\documentclass[11pt,a4paper,titlepage,fleqn,]{report}
\usepackage[centertags]{amsmath}
\usepackage{amsfonts}
\usepackage{amssymb}
\usepackage{amsthm}
\usepackage{hyperref}
\usepackage{newlfont}
\usepackage{verbatim}
\newtheorem{lem}{Lemma}[chapter]
\newtheorem{cor}[lem]{Corollary}
\newtheorem{theo}[lem]{Theorem}
\newtheorem{defi}[lem]{Defination}
\begin{document}
\title{ON THE EIGENVALUE OF $p(x)$-LAPLACE EQUATION}
\author{Yushan Jiang, Yongqiang Fu}
\maketitle
\begin{abstract}

The main purpose of this paper is to show that there exists a
positive number
    $\lambda_{1}$,
the first eigenvalue, such that some $p(x)$-Laplace equation admits
a solution if
    $\lambda=\lambda_{1}$
and that
    $\lambda_{1}$
is simple, i.e., with respect to \textit{the first eigenvalue}
solutions, which are not equal to zero a. e., of the $p(x)$-Laplace
equation forms an one dimensional subset. Furthermore, by developing
Moser method we obtained some results concerning H\"{o}lder
continuity and bounded properties of the solutions. Our works are
done in the setting of the Generalized-Sobolev Space. There are many
perfect results about $p$-Laplace equations, but about
$p(x)$-Laplace equation there are few results. The main reason is
that a lot of methods which are very useful in dealing with
$p$-Laplace equations are no longer valid for $p(x)$-Laplace
equations. In this paper, many results are obtained by imposing some
conditions on $p(x)$.

Stimulated by the development of the study of elastic mechanics,
interest in variational problems and differential equations has
grown in recent decades, while Laplace equations with nonstandard
growth conditions share a part. The equation discussed in this paper
is derived from the elastic mechanics.\newline
\end{abstract}
\tableofcontents
\chapter{Introduction}
\section{Background}
In recent years, there has been increasing interest toward variable exponent Lebesgue
and Sobolev spaces.It is clear that we cannot simply replace
    $p$
by
    $p(x)$
in the usual definition of the norm in
    $L^{p}$.
However, the Lebesgue spaces can be considered as particular cases of the Orlicz
spaces belong to a larger family of so called modula spaces. This approach enables to
define corresponding counterparts of the Luxemburg and Orlicz norms in
    $L^{p(x)}$.
The present line of investigation toward variable exponent Lebesgue and Sobolev
spaces goes back to a paper by O.Kov\'{a}$\breve{c}$ik and
J.R\'{a}kosn\'{\i}k
    \cite{O.Ko}
from 1991.After this paper not much happened till the late 1990's. At this point the
subject seems to have been rediscovered by several researchers independently:S.Samko
    \cite{S.Sa,S.Sam},
working based on earlier Russian work(I.Sharapudinov and V.Zhikov
    \cite{V.V.Z}),
X.Fan and collaborators drawing inspiration from the study of differential equations
    \cite{X.Fa.Embed,X.Fa.Quasi,X.Fa.Space,X.Fa.Compa}.
The last couple of years have seen the integration of the separate lines of
investigation, but much still remains to be done.

The main incentive for many of the investigators of variable exponent spaces is
relaxing the coercivity conditions assumed for the solutions of a differential
equation or the corresponding variational integral.while Laplace equations with
nonstandard growth conditions share a part.One such application has been investigated
in greater detail, electro-rheological fluids. These fluids have the interesting
property that their viscosity depends on the electric field in the fluid. For some
technical applications the mathematical theory was presented by some
investigators:M.R$\check{u}\check{z}$i$\check{c}$ka, E.Acerbi and G.Mingione
    \cite{E.Ace.01,E.Ace.02,E.Ace.2,M.Ruz}.

The remain part of this part strive to give a little more detailed
an account of the mathematics of variable exponent spaces and in the
next section we present the theory of Generalized Lebesgue spaces
and that of Generalized Sobolev spaces. In the final section or the
main section we present a generalization of the eigenvalue problem
on some $p(x)$-Laplace equation by Mitsuharu \^{O}tani and
Toshiak\'{\i}  Teshima
    \cite{T.Teshima}.
\section{Overview of Differential Operator with Nonstandard Exponent Growth}
\subsubsection{The Harnack Inequality}
In a bounded domain
    $D$
of Euclidean space
    $\mathbb{R}^{n},n\geq2$,
Yu.A.Alkhutov\cite{Yu.A.A} proved the Harnack inequality and an interior a priori
estimate for the H\"{o}lder norm of solutions about the equation as following:
    \begin{equation}
    \sum^{n}_{i=1}
        \frac{\partial}{\partial x_{i}}
    \left(
        |\nabla u|^{p(x)-2}
        \frac{\partial u}{\partial x_{i}}
    \right)
    =0
    \end{equation}
where
    $p(x)$
is a measurable function in
    $D$
and
    $1<p_{1}\leq p(x)\leq p_{2}<\infty$,
the domain
    $D$
is divided by a part of a Lipschitz surface into two subdomains, in each of which
    $p(x)$
is constant.
\subsubsection{The Hardy-Littlewood Maximal Operator.}

Assume that
    $1<p^{-}\leq p^{+}<\infty$
and there exists a constant
    $C>0$
such that
    \begin{displaymath}
        |p(x)-p(y)|\leq\frac{C}{-\log|x-y|}
    \end{displaymath}
for every
    $x,y\in \mathbb{R}^{n}, |x-y|\leq\frac{1}{2}$
and
    \begin{displaymath}
        |p(x)-p(y)|\leq\frac{C}{-\log(e+|y|)}
    \end{displaymath}
for every
    $x,y\in \mathbb{R}^{n},|y|\geq |x|$.
Under these assumption on
    $p(x)$,
Cruz-Uribe, Fiorenza and Neugebauer
    \cite{Cru and Fior}
proved that the Hardy-Littlewood maximal operator is bounded from
    $L^{p(\cdot)}(\mathbb{R}^{n})$
to itself. This was an improvement of earlier work by Diening
    \cite{Die.Max,Die.Embed,Die.Ruz.Ca,Die.Ruz.In}
and Nekvinda
    \cite{A.Nekv}.
Maximal operator have also been studied in weighted
    $L^{p(\cdot)}$
spaces by Kokilashvili and Samko
    \cite{V.Koki,S.Sam,S.Sa}.
\subsubsection{Strong Maximum Principle of $p(x)$-Laplace Equation}
If
    $\Omega\subset \mathbb{R}^{N} (N\geq2)$
be an open set,
    $p(x)\in C^{1}(\overline{\Omega})$,
and
    $p(x)>1(x\in \overline{\Omega}), q(x)\in C^{0}(\overline{\Omega})$,
and
    $p(x)\leq q(x)\leq p^{\ast}(x)$
    ($p^{\ast}=\frac{Np(x)}{N-p(x)}$
for
    $p(x)<N$;
    $p^{\ast}=\infty$
for
    $p(x)\geq N$),
    $d(x)\in L^{\infty}(\Omega),d(x)\geq 0\ a.e..$
Fan.X and Zhao.Y
    \cite{X.Fa.Strong}
given a strong maximum principle for super-solutions of the $p(x)$-Laplace equations
    \begin{equation}\label{eq:fan}
        -div(|\nabla u|^{p(x)-2}\nabla u)+d(x)|u|^{q(x)-2}u
        =0
    \end{equation}
Fan.X and Zhao.Y proved that the nonnegative weak upper solution of
    (\ref{eq:fan})
    $u$
satisfies
    $u\geq c, x\in K\ a.e.$,
for any given nonempty compact subset
    $K\subset\Omega$,
where
    $c>0$
is a constant. Furthermore if
    $u\in C^{1}(\Omega \cup {x_{1}}),u(x_{1})=0, x_{1}\in \partial \Omega$
and
    $u$
satisfies the inner sphere conditions, then
    $\frac{\partial u(x_{1})}{\partial\gamma}>0$,
where
    $\gamma$
is the unit inner normal vector of
    $\partial\Omega$
at
    $x_{1}$.
\subsubsection{Existence of Solutions for Elliptic Systems with Nonuniform Growth}
Let
    $\Omega\subset\mathbb{R}^{n}$
be a bounded Lipschitz domain. For the following systems:
    \begin{equation}\label{eq:fu0}
        \frac{\partial A^{i}_{\alpha}}{\partial A^{\alpha}}
        (x,u(x),Du(x))
        =
        B^{i}(x,u(x),Du(x)),
        x\in\Omega, i=1,\ldots,N
    \end{equation}
    \begin{equation}\label{eq:fu1}
        u^{i}(x)=0,\ x\in \partial\Omega, \qquad i=1,\ldots,N
    \end{equation}
Fu.Yongqiang proved that if the coeifficients of
    (\ref{eq:fu0})
satisfy
    $A^{i}_{\alpha}:
    \Omega\times\mathbb{R}^{N}\times M^{N\times n}
    \rightarrow \mathbb{R},
    B^{i}:
    \Omega\times\mathbb{R}^{N}\times
    M^{N\times n}\rightarrow \mathbb{R},
    i=1,\ldots,N,\alpha=1,\ldots,n,$
are Carath\'{e}odory functions.

    $
    |A(x,s,\xi)|
    \leq
    C_{1}|\xi|^{p(x)-1}+C_{2}|s|^{p(x)-1}+G(x)
    $,
 where
    $G\in L^{p'(x)}(\Omega)$,
    $C_{1},C_{2}\geq 0$
 and
    $C_{2}$ small.

    $
    |B(x,s,\xi)|
    \leq
    C'_{1}|\xi|^{p(x)-1}+C'_{2}|s|^{p(x)-1}+\overline{G}(x)
    $,
 where
    $\overline{G}\in L^{p'(x)}(\Omega)$,
    $C'_{1},C'_{2}\geq 0$
 and
    $C'_{2}$
 small.

    $A^{i}_{\alpha}(x,s,\xi)\xi^{i}_{\alpha}
    \geq
    \lambda_{0}|\xi|^{p(x)}-C|s|^{p(x)}+h(x)$,
 where
    $\lambda_{0}>0,C\geq 0$
small and
    $h\in L^{1}(\Omega)$.

 For almost every
    $x_{0}\in \Omega, s_{0}\in\mathbb{R}^{N},$
 the mapping
    $\xi\mapsto A(x_{0},s_{0},\xi)$
 satisfies
    \begin{displaymath}
    \int_{G}A^{i}_{\alpha}(x_{0},s_{0},\xi_{0}+Dz(x))z^{i}_{\alpha}(x)dx
    \geq
    \nu\int_{G}|Dz(x)|^{p(x)}dx
    \end{displaymath}
for each
    $\xi_{0}\in M^{N\times n}$,
    $G\subset \mathbb{R}^{n}$,
    $z\in C^{1}_{0}(G,\mathbb{R}^{N})$
where
    $\nu >0$
and
    $(Du(x))^{i}_{\alpha}
    =
    \partial u^{i}(x)/\partial{x^{\alpha}}
    =
    u^{i}_{\alpha}(x)$.
    $p:\Omega \rightarrow [1,\infty]$
is a measurable function and $p'$ is its conjugate function. Then the Dirichlet
problem
    (\ref{eq:fu0}),(\ref{eq:fu1})
has at least one weak solution in
    $W^{1,p(\cdot)}_{0}(\Omega,\mathbb{R}^{N})$,
that is to say, there exists at least one
    $u\in W^{1,p(\cdot)}_{0}(\Omega,\mathbb{R}^{N})$
satisfying
    \begin{equation}
    \int_{\Omega}
    [
    A^{i}_{\alpha}(x,u,Du)z^{i}_{\alpha}(x)
    +
    B^{i}(x,u,Du)z^{i}(x)
    ]dx
    =0
    \end{equation}
for all
    $x\in W^{1,p(\cdot)}_{0}(\Omega,\mathbb{R}^{N})$.
This generalizes the result of Acerbi and Fusco
    \cite{E.Ace and N.Fusco}.
\subsubsection{H\"{o}lder Continuity of Minimizers of Functionals
               with Variable Growth Exponent}
Let
    $a\in W^{1,s}(\Omega)(s>n), r$
be two nonnegative measurable functions such that
    $1<p_{0}\leq a(x)\leq q_{0}\leq p^{\ast}_{0}$,
    $0\leq r(x)\leq r\leq p^{\ast}_{0}$
and let
    $f:\Omega\times\mathbb{R}\times\mathbb{R}^{n}$
be a Carath\'{e}odory function satisfying the growth assumptions
    $c_{1}(|\xi|^{a(x)}-|u|^{r(x)}-1)
    \leq
    f(x,u,\xi)
    \leq
    c_{2}(|\xi|^{a(x)}+|u|^{r(x)}+1)$
and let
    $u$
be a quasiminimizer of the functional as following
    \begin{equation}
        \mathcal{F}(u)
        =
        \int_{\Omega}f(x,u,Du)dx
    \end{equation}
Valeria Chiad\`{o}Piat and Alessandra Coscia\cite{Vale}proved the
the locally H\"{o}lder continuous of $u$ in $\Omega$.
\subsubsection{H\"{o}lder continuity of $p(x)$-Laplace equation}
Let
    $\Omega$
be a open set in
    $\mathbb{R}^{N}$,
for the following equation:
\begin{equation}\label{eq:fan2}
    -div\left(\lambda+|\nabla u|^{2}\right)^{\frac{p(x)-2}{2}}\nabla u
    =
    F(x,u),
    x\in\Omega\subset\mathbb{R}^{N}
\end{equation}
where
    $\lambda
    \geq
    0,
    F
    \in
    C^{0}(\Omega\times\mathbb{R})$
satisfies
    $|F(x,u)|
    \leq
    c_{1}+c_{2}|u|^{q(x)},
    \forall
    (x,u)\in\Omega\times\mathbb{R}$
where
    $1<q(x)<p^{\ast}(x)
    (p^{\ast}
    =
    \frac{Np(x)}{N-p(x)}
    \ for\  p(x)<N,
    p^{\ast}
    =
    +\infty
    \ for
    p(x)\geq N)$,
    $p\in C^{1}(\Omega),
    p(x)>1
    (\forall x\in \Omega)$.
X.Fan and Zhao Dun
    \cite{X.Fan Regu}
proved the local
    $C^{1,\alpha}$
regularity $u$, that is to say, the weak solution of
    (\ref{eq:fan2})
satisfies
    $u\in C^{1,\alpha}_{Loc}(\Omega)$.
\subsubsection{On the Positive Solution of $p(x)$-Laplace Equation}
Let
    $\Omega$
be a bounded domain in
    $\mathbb{R}^{N}(N>1)$,
for the following equation:
    \begin{equation}\label{eq:exi}
    \left\{ \begin{array}{ll}
    -div(|\nabla u|^{p(x)-2}\nabla u)
    =\lambda|u|^{\alpha(x)-2}u+|u|^{\beta(x)-2}u&\ x \in\Omega \\
    u(x)=0 &\ x \in \partial\Omega
    \end{array} \right.
    \end{equation}
where
    $p,\alpha,\beta$
are the continuous functions on
    $\overline{\Omega}$,
and
    $p(x)<N,\lambda >0$.
X.Fan\cite{X.Fan Exi} proved that if
    1)
    $p(x):\overline{\Omega}\rightarrow \mathbb{R}$
is Lipschitz continuous and
    $p_{-}>1$;
    2)
    $1<\alpha_{-}\leq \alpha^{+}<p_{-}\leq p^{+}<\beta_{-},
    \beta(x)\leq p^{\ast}(x)$.
then (\ref{eq:exi}) has at least two positive solutions for small
    $\lambda$.
\subsubsection{Dirichlet Boundary Value Problem}
Consider a differential operator $A$ of order $2k$ in the divergence form
    \begin{equation}
        Au(x)
        =
        \sum_{|\alpha|\leq k}(-1)^{|\alpha|}D^{\alpha}a_{\alpha}(x,\delta_{k}u(x))
    \end{equation}
where the functions
    $a_{\alpha}(x,\delta_{k}u(x))\in CAR(\Omega, m), $
    $m=\sharp \{\alpha\in \mathbb{N}^{N}_{0}:|\alpha|\leq k\}$,
fulfill the growth condition
    $|a_{\alpha}(x, \xi)|
    \leq
    g(x)
    +
    c\sum_{|\alpha|\leq k}|\xi_{\alpha}|^{p(x)-1}$
with
    $g\in C^{ps(x)}(\Omega)$
and
    $c>0$.
Let
    $Q$
be a Banach space of functions on $\Omega$ equipped with a norm
    $\|\cdot\|_{Q}$
and such that
    $C^{\infty}_{0}(\Omega)$
is dense in $Q$ and moreover,
    $W^{k,p(x)}_{0}(\Omega)\circlearrowleft Q$
A function
    $u\in W^{k,p(x)}(\Omega)$
is a weak solution to the Dirichlet boundary value problem
    $(A,u_{0},f)$
for the equation
    $Au=f$
with the boundary condition given by
    $u_{0}$,
if
    $u-u_{0}\in W^{k,p(x)}_{0}(\Omega)$
and if the identify
    \begin{equation}
        \sum_{|\alpha|\leq k}
        \int_{\Omega}
            a_{\alpha}(x,\delta_{k}u(x))D^{\alpha}\upsilon(x)
        dx
        =
        \langle f,\upsilon\rangle
    \end{equation}
holds for every
    $\upsilon\in W^{k,p(x)}_{0}(\Omega)$.
O.Kov\'{a}$\breve{c}$ik and J.R\'{a}kosn\'{\i}k
    \cite{O.Ko}
proved that if
    $p(x)\in P(\Omega)$
satisfy
    \begin{displaymath}
        1<ess\inf_{\Omega}p(x)\leq ess\sup_{\Omega}p(x)<\infty
    \end{displaymath}
the functions $a_{\alpha}$ satisfy
\begin{equation}\label{ineq:o.kav1}
    \sum_{|\alpha|\leq k}
        \left[
            a_{\alpha}(x,\xi)
            -
            a_{\alpha}(x,\eta)
        \right]
    (\xi_{\alpha}-\eta_{\alpha})
    \geq
    0,
\end{equation}
    \begin{equation}\nonumber
    \sum_{|\alpha|\leq k}
        a_{\alpha}(x,\xi)\xi_{\alpha}
    \geq
    c_{1}\sum_{|\alpha|\leq k}
        |\xi_{\alpha}|^{p(x)}-c_{2}
    \end{equation}
for every
    $\xi,\eta \in \mathbb{R}^{m}$
and for a.e.
    $x\in\Omega$
with some constants
    $c_{1},c_{2}>0$.
Then the boundary value problem
    $(A,u_{0},f)$
has at least one weak solution
    $u\in W^{k,p(x)}(\Omega)$.
If moreover, the inequality
    (\ref{ineq:o.kav1})
is strict for
    $\xi\neq \eta$
then the solution is unique.
\subsubsection{Existence of Positive Solution on $p(x)$-Laplace Equation}
For the equations as following:
    \begin{equation}\label{eq:fan1}
    \left\{
        \begin{array}{ll}
        -div(|\nabla u|^{p(x)-2}\nabla u)
        =
        a(x)|u|^{\beta(x)}u&\ x\in\Omega \\
        u(x)
        =
        0 &\ x \in \partial\Omega
        \end{array}
    \right.
    \end{equation}
X.Fan
    \cite{X.Fa.Exi2}
proved the existence of positive solution for
    (\ref{eq:fan1})
with
    $\beta^{+}<p_{-}$
or
$\beta_{-}>p^{+}$.
\section{Origin of Problem}
Our study comes from the article
    \cite{T.Teshima}
written by M.$\hat{o}$tani and T.Teshima. in which they study the eigenvalue of the
equation as following:
    \begin{displaymath}
        \left\{
            \begin{array}{ll}
                -\Delta_{p}u(x)+a(x)|u(x)|^{p-2}u(x)
                =
                \lambda b(x)|u(x)|^{p-2}u(x)&\ x \in\Omega \\
                u(x)
                =
                0 &\ x \in \partial\Omega
            \end{array}
        \right.
    \end{displaymath}
where
    $\Delta_{p}u(x)
    =
    Div(|\nabla u(x)|^{p-2}\nabla u(x)),
    \lambda>0$.

M.$\hat{o}$tani and T.Teshima. proved that eigenvalue problem above
has a nontrivial nonnegative solution $u$ if and only if
    $\lambda=\lambda_{1}$
and
    $J_{\lambda_{1}}:
    =
    A(u)-\lambda_{1}B(u)
    =
    0$.
Furthermore, the set of all solutions consists of
    $tu_{1};t\in \mathbb{R}^{1}$
where
    $u_{1}$
is a solution and
    $u_{1}\in C^{1,\theta}(\overline{\Omega})$
for some
    $\theta\in (0,1)$.
But if
    $p$
is a function of
    $x\in\Omega$.
it's a more difficult situation. It's clear that we can not simply replace $p$ by
$p(x)$ in the equation about
    \cite{X.Fa.Compa}.
However, we can extend the definition of $p$-Laplace operator by
    \begin{equation}\label{def:p(x)-Lap.Oper}
        \Delta_{p(x)}u(x)
        =
        Div(p(x)|\nabla u(x)|^{p(x)-2}\nabla u(x))
    \end{equation}
which is named $p(x)$-Laplace Operator. Our work tries to give some results about the
eigenvalue problem of $p(x)$-Laplace equation as following:
\begin{equation}\label{eq:jys}
     -\Delta_{p(x)}u(x)+a(x)|u(x)|^{p(x)-2}u(x)
     =\lambda b(x)|u(x)|^{p(x)-2}u(x)
\end{equation}
\chapter{Generalized Sobolev Space}
\section{Conceptions and Properties}
For a set
    $\Omega\in \mathbb{R}^{N}$
with
    $|\Omega|>0$,
we define the family of all measurable functions
    $p:\Omega\rightarrow[1,\infty]$
by
    $\mathcal{P}(\Omega)$.
    we put
    $\Omega_{1}
    =
    \{
        x\in\Omega:p(x)=1
    \}$,
    $\Omega_{\infty}
    =
    \{
        x\in\Omega:p(x)=\infty
    \}$,
    $\Omega_{0}
    =
    \Omega\backslash(\Omega_{0}\cup\Omega_{\infty})$;
also, we define
    $p^{+}
    =
    \textrm{ess}\sup_{x\in\Omega_{0}}p(x)$
and
    $p^{-}
    =
    \textrm{ess}\inf_{x\in\Omega_{0}}p(x)$.
We define the variable exponent Lebesgue space
    $L^{p(\cdot)}(\Omega)$
to consist of all measurable functions
    $u:\mathbb{R}^{n}\rightarrow\mathbb{R}$
such that
    \begin{equation}
        \varrho_{p(\cdot)}(\lambda u)
        =
        \int_{\mathbb{R}^{n}}
            |\lambda u(x)|^{p(x)}
        dx
        +
        \textrm{ess}\sup_{\Omega_{\infty}}|u(x)|
        <
        \infty
    \end{equation}
for some $\lambda >0$.
    The function
    $\varrho_{p(\cdot)}:
    L^{p(\cdot)}(\Omega)\rightarrow[0,\infty]$
is called the modular of the space
    $L^{p(\cdot)}(\Omega)$.
we define a norm, the so-called Luxemburg norm, on this space by the formula
    \begin{equation}
        \|u\|_{p(\cdot)}
        =
        \inf
        \{
            \lambda>0:\varrho_{p(\cdot)}(u/\lambda)
            \leq
            1
        \}
    \end{equation}
If $p$ is a constant function, then the variable exponent Lebesgue spaces coincides
with the classical Lebesgue spaces and so the notation can give rise to no confusion.
The varialbe exponent Sobolev space
    $W^{k,p(x)}(\Omega)$
is the subspace of functions
    $L^{p(\cdot)}(\Omega)$
whose distributional gradient exists almost everywhere and lies in
    $L^{p(\cdot)}(\Omega)$.
The norm of
    $W^{k,p(x)}(\Omega)$
defined by
\begin{equation}\label{def:sb.norm}
    \|u\|_{k,p}
    =
    \sum_{|\alpha|\leq k}\|D^{\alpha}u\|_{p(\cdot)}
\end{equation}
By
    $W^{k,p(x)}_{0}(\Omega)$
we denote the subspace of
    $W^{k,p(x)}(\Omega)$
which is the closure of
    $C^{\infty}_{0}$
with respect to the norm
    (\ref{def:sb.norm}).
\newline
    \textbf{Basic properties.}
Variable exponent Lebesgue spaces resemble classical Lebesgue spaces in many respects
    \cite{O.Ko}
--they are Banach spaces,the H\"{o}lder inequality holds,they are reflexive if and
only if
    $1<p^{-}\leq p^{+}<\infty$
and continuous functions are dense if
    $p^{+}<\infty$;
    The inclusion between Lebesgue spaces also generalizes naturally:if
    $0<|\Omega|<\infty$
and $p,q$ are variable exponent so that
    $p(x)\leq q(x)$
almost everywhere in $\Omega$ then there exists an imbedding
    $L^{q(\cdot)}(\Omega)\hookrightarrow L^{p(\cdot)}(\Omega)$
whose norm does not exceed
    $|\Omega|+1$;
    If
    $p^{+}<\infty$
and
    $(f_{i})$
is a sequence of functions in
    $L^{q(\cdot)}(\Omega)$,
then
    $\|f_{i}\|_{p(\cdot)}\rightarrow 0$
if and only if
    $\varrho_{p(\cdot)}(f_{i})\rightarrow 0$;
    The spaces
    $W^{k,p(x)}(\Omega)$ and $W^{k,p(x)}_{0}(\Omega)$
are Banach spaces,which are separable and reflexive if
    $1<p^{-}\leq p^{+}<\infty$;
If
    $q(x)\leq p(x)$ for a.e. $x\in \Omega$
then
    $W^{k,p(x)}(\Omega)\circlearrowleft W^{k,q(x)}(\Omega)$.
\section{Sobolev Embedding Inequalities}
As we know,in dealing with some partial differential equation
problems Sobolev embedding inequality is useful for us. Many good
results are derived from these inequalities.
    \begin{theo}[\cite{X.Fa.Embed}]Assum $\Omega$ be a open domain in
            $\mathbb{R}^{N}$,
        with cone property let
            $p(x)\in \mathcal{P}(\overline\Omega)$
        be a Lipschitz continuous function, if
            $q(x)\in \mathcal{P}(\Omega)$
        satisfy
            $p(x)\leq q(x)\leq P^{\ast}(x):
            =\frac{Np(x)}{N-kp(x)} a.e. x\in\overline\Omega$
        then there exists a continuous embedding
            $W^{k,p(x)}(\Omega)\hookrightarrow L^{q(x)}(\Omega)$.
    \end{theo}
    \begin{theo}[\cite{X.Fa.Embed}]
        \label{emb:2} Assum $\Omega$ be a open domain in
            $\mathbb{R}^{N}$,
        with cone property;
        If
            $p(x):\overline\Omega\rightarrow \mathbb{R}$
        is uniform continuous and satisfy
            $1<p^{-}\leq p^{+}<\frac{N}{k}$
        then for any measurable function $q(x)$ defined in $\Omega$ with
            $p(x)\leq q(x), a. e. x\in \overline\Omega$
        and
            $\textrm{ess}\inf_{x\in\overline\Omega}(p^{\ast}(x)-q(x))>0$
        there is a continuous embedding
            $W^{k,p(x)}(\Omega)\hookrightarrow L^{q(x)}(\Omega)$.
    \end{theo}
    \begin{theo}[\cite{X.Fa.Embed}]
        Assume that $\Omega$ be a open domain in
            $\mathbb{R}^{N}$
        with cone property,If
            $\Omega$ is bounded,
            $p(x)\in C(\overline\Omega)$
        and
            $q(x)$ is the same as in theorem \ref{emb:2},
        then there is a continuous compact embedding
            $W^{k,p(x)}(\Omega)\circlearrowleft L^{q(x)}(\Omega)$.
    \end{theo}
    \begin{theo}[\cite{X.Fa.Compa}]
        Suppose that
            $p:\mathbb{R}^{N}\rightarrow \mathbb{R}$
        is a uniformly continuous and radically symmetric function
        satisfying
            $1<p^{-}\leq p^{+}<N$
        then, for any measurable function
            $\alpha :\mathbb{R}^{N}\rightarrow \mathbb{R}$
        with
            $p(x)\ll \alpha(x)\ll p^{\ast}(x),\forall x\in\mathbb{R}^{N}$,
        there be a compact  imbedding
            $W^{1,p(x)}_{r}(\mathbb{R}^{N})\circlearrowleft
            L^{\alpha(x)}(\mathbb{R}^{N})$.
        where
            $W^{1,p(x)}_{r}(\mathbb{R}^{N})
            :=\{u\in W^{1,p(x)}(\mathbb{R}^{N}):u\textrm{is radially
            symmetric.}\}$.
    \end{theo}
    \begin{theo}[\cite{X.Fa.Compa}]
        If
            $p:\mathbb{R}^{N}\rightarrow \mathbb{R}$
        is a uniformly continous and satisfies
            $1<p^{-}\leq p^{+}<N$
        then for any measurable function $\alpha$ with
            $\overline p\ll\alpha\ll \overline p^{\ast}, x\in \mathbb{R}^{N}$
        we have the compact imbedding
            $W^{1,p(x)}_{r}(\mathbb{R}^{N})\circlearrowleft
            L^{\alpha(x)}(\mathbb{R}^{N})$.
        where
            $\overline p(x)=\int_{G}p(g(x))d\mu(g),\forall x\in
            \mathbb{R}^{N}$,
        $G=O(N)$ be the orthogonal group on $\mathbb{R}^{N}$, $\mu$
        be a Haar measure on the compact group $G$, and $\mu (G)=1$.
    \end{theo}
    \begin{theo}[\cite{X.Fa.Compa}]
         Let $G$ be a subgroup of $O(N)$ and
        $\Omega$ be a invariant open subset in $\mathbb{R}^{N}$
        compatible with $G$,
            $p:\mathbb{R}^{N}\rightarrow \mathbb{R}$
        be G-invariant and uniformly continuous such that
            $1<p^{-}\leq p^{+}<N$
        holds. Then ,for any measurable function $\alpha$ with
            $ p\ll\alpha\ll  p^{\ast}, x\in \Omega$,
        we have the compact imbedding
            $W^{1,p(x)}_{0,G}(\Omega)\circlearrowleft
            L^{\alpha(x)}(\Omega)$.
        where
            $W^{1,p(x)}_{0,G}(\Omega):=
            \{u\in W^{1,p(x)}_{0}(\Omega):u \textrm{ is
            G-invariant}\}$.
    \end{theo}
\section{Notations and Preliminary Results}
Let $\Omega$ be a bounded open subset of
    $\mathbb{R}^{n},(n\geq 2)$.
we use the standard notation for the Generalize Lebesgue and
Generalize Sobolev spaces
    $L^{p(x)}(\Omega)$
and
    $W^{k,p(x)}(\Omega)$;
in particular we will denote by
    $\|\cdot\|_{p}$
and
    $\|\cdot\|_{k,p}$
the corresponding norms. The Lebesgue measure of a measurable set
    $A\subset\mathbb{R}^{n}$
will be denoted by $|A|$,whereas a ball of radius $R$ will be denoted
by $B_{R}$ and all balls mentioned in a single proposition will
always be assumed to be concentric;moreover if
    $u:\Omega\rightarrow\mathbb{R},k\in\mathbb{R}$
and
    $B_{R}$
is a ball strictly contianed in $\Omega$, we set
\begin{eqnarray*}
  A(k,R) &=& \{x\in B_{R}:u(x)>k\} \\
  M(u,R) &=& \sup_{B_{R}}u \\
  m(u,R) &=& \inf_{B_{R}}u \\
  osc(u,R) &=& M(u,R)-m(u,R)
\end{eqnarray*}
Finally ,for all $k\in\mathbb{R}$ and $R>0$ we set
\begin{eqnarray*}
  \oint_{B_{R}}udx &=& \frac{1}{|B_{R}|}\int_{B_{R}}udx \\
  \oint_{A(k,R)}udx &=& \frac{1}{|B_{R}|}\int_{A(k,R)}udx
\end{eqnarray*}
we will denote by the same letter $C$ (or $C(\cdots)$ to stress the
dependence on some arguments) several constants, whose value may
change from line to line.

In the proof of the article we shall use the following Lemma which
can be found in \cite{Chen.y.z},\cite{Yongq},\cite{Vale}.
\begin{lem}[Moser iteration inequality]\label{ineq:Moser}
Let $\{x_{i}\}_{i}$ be a sequence of positive real numbers such that
    $x_{0}\leq C^{-1/\beta B^{-1/\beta^{2}}}$,
    $x_{i+1}\leq CB^{i}x_{i}^{1+\beta}$
where
    $\beta>0,C>0,B>1.$
then
    $x_{i}\rightarrow 0$ as $i\rightarrow +\infty$.
\end{lem}
\begin{lem}[Poicar\'{e} inequality\cite{Chen.y.z}]
Let $u$ be a function in $W^{1,p}_{0}(\mathbb{R}^{n})$ then
    \begin{equation}\label{ineq:pocar}
        \int_{\mathbb{R}^{n}}|u|^{p^{\ast}}dx
        \leq \prod^{n}_{i=1}\left\{\left(\int_{\mathbb{R}^{n}}
        |u|^{p}dx\right )^\frac{1}{p} \right\}^{\frac{p^{\ast}}{n}}
    \end{equation}
holds for $1\leq p\leq n$,where $p^{\ast}=\frac{np}{n-p}$.
\end{lem}
\begin{lem}[Sobolev-Poicar\'{e} inequality\cite{Yongq}]
    \label{ineq:s-p}
for any given bounded domain $\Omega$,if
    $p(x)\in L^{\infty}(\Omega), u(x)\in W^{1,p(x)}_{0}(\Omega)$
then
    \begin{equation}
        \int_{\Omega}|u(x)|^{p(x)}dx
        \leq C\int_{\Omega}|\nabla u(x)|^{p(x)}dx
    \end{equation}
where $C$ is a constant depended on $\Omega$.
\end{lem}
\begin{lem}[Sobolev-Poicar\'{e} inequality with variable exponent\cite{Vale}]
\label{lemma:ineq:s-p with nonstandard}
Assume that
    $a\in W^{1,s}(\Omega),s>n$
satisfies
    $1<p_{0}<a(x)<q_{0}\leq p^{\ast}_{0}$;
then for every $M>0$ there exists a positive radius
    $R_{1}=R_{1}(M,s,n,\|a\|_{1,s})$
such that for every
    $\gamma>\frac{1}{n}-\frac{1}{s}$
there exist two positive constants
    $\chi=\chi(n,p_{0},s,\gamma,\|a\|_{1,s})$
and
    $C=C(n,p_{0},q_{0})$
for which
\begin{equation}\label{ineq:s-p varia}
    \left(\oint_{B_{R}}\left|\frac{u}{R}\right|^{a(x)\frac{n}{n-1}dx}\right)^{\frac{n-1}{n}}
    \leq c\oint_{B_{R}}|Du|^{a(x)}dx+\chi|\{x\in B_{R}:|u|>0\}|^{\gamma}
\end{equation}
holds for every $B_{R}\subset\Omega$ with $0<R<R_{1}$,and every $u\in
W^{1,p_{0}}(B_{R})$ such that
    $\int_{B_{R}}|Du|^{a(x)}dx<+\infty, \sup_{B_{R}}\leq M\
    \ u=0\ \textrm{on}\ \partial B_{R}$.
\end{lem}
\chapter{On the First Eigenvalue Problem of $p(x)$-Laplace Equation}

\section{Introduction}
    Suppose
        $\Omega$
    is a bounded domain in
        $R^{n}$
    with a smooth boundary
        $\partial\Omega$.
    For certain given
        $p(x)\in P(\Omega)$,
    where
        $1<p\leq p(x)\leq p^{*}< +\infty (p^{*}=\frac{np}{(n-p)})$
    and
        $p(x)$
    is continuous in $\Omega$ and the partial differential of
        $p(x)$, $\frac{\partial p}{\partial x_{i}}$
    is bounded a.e. in $\Omega$.
Thinking the eigenvalue problem of $p(x)$-Laplace equation
    $(E)_{\lambda}$
as following:
    \begin{displaymath}
        \left\{
            \begin{array}{ll}
                -\Delta_{p(x)}u(x)+a(x)|u(x)|^{p(x)-2}u(x)=
                \lambda b(x)|u(x)|^{p(x)-2}u(x)&\ x \in\Omega \\
                u(x)=0 &\ x \in \partial\Omega
            \end{array}
        \right.
    \end{displaymath}
where
    \begin{eqnarray}
        &&
        -\Delta_{p(x)}u(x)
        =
        -Div(p(x)|\nabla u(x)|^{p(x)-2}\nabla u(x)) \\
        &&
        a(x)\in L^{\infty}_{+}(\Omega)
        =
        \{f\in L^{\infty}(\Omega)|f(x)\geq 0 a.e. x\in\Omega\} \\
        &&
        b(x)\in L^{\infty}_{0}(\Omega)
        =
        \{f\in L^{\infty}(\Omega)|f^{+}(x)
        =:max\{f(x),0\}\neq 0 \}
    \end{eqnarray}
We say that
    $u(x)$
is the solution of the eigenvalue problem
    $(E)_{\lambda}$
if
    $u(x)\in W^{1,p(x)}_{0}(\Omega)$
satisfies the equation in the general sense, that is to say, for any given function
    $h(x)\in C^{\infty}_{0}(\Omega)$
there stands
    \begin{eqnarray}\label{def:weak solution}
        &&
        \int_{\Omega}
        -Div(p(x)|\nabla u(x)|^{p(x)-2}\nabla u(x))h(x)
        dx \nonumber \\
        &&+\int_{\Omega}a(x)|u(x)|^{p(x)-2}u(x)h(x)dx \nonumber \\
        &&=\lambda\int_{\Omega}b(x)|u(x)|^{p(x)-2}u(x)h(x)dx
    \end{eqnarray}
In the generalized sobolev space, we study the eigenvalue of the
equation.By Moser iteration, we obtain some properties about the
solution of eigenvalue problem $(E)_{\lambda}$:boundary,H\"{o}lder
continouty,and so on.

\section{On the Local Boundary of Eigenvalue Problem}

\subsection{Gradient Estimate of Solution}

    \begin{lem}[Caccioppoli Inequality on  Solution]
        Suppose
            $u(x)$
        is the solution of
            $(E)_{\lambda}$, $p(x)$
        satisfies
            $1 < p \leq p(x) \leq p <  p^{*}<+\infty$
        in a certain open set
            $\omega$
        included in
            $\Omega$,
        then for any  given spherical neighborhood
            $B_{R}$
        in
            $\omega$
        and every
            $0<s<t<R<1$ , $k>0$
        there stands
            \begin{equation}\label{ineq:cacci}
                \int_{A(k,s)}
                    \left|
                        \nabla u(x)
                    \right|^{p}dx
                \leq C
                \left(
                    \int_{A(k,t)}
                        \left|
                            \frac{u(x)-k}{t-s}
                        \right|^{q}dx+(1+k^{q})|A(k,R)|
                \right)
            \end{equation}
        where
            $C=C(p,q,\lambda)$
    \end{lem}

\textbf{Proof}:
    Let
        $\zeta(x)$
    be a cut-off function between
        $B_{s}$ and $B_{t}$
    with
        $|\nabla\zeta(x)|\leq\frac{2}{t-s}$,
    by using a test function
        $\varphi=-\zeta(u-k)^{+}$
    in equation (\ref{eq:e}), one arrives at
        \begin{equation*}
                \int_{A(k,s)}
                    \left|
                        \nabla u(x)
                    \right|
                ^{p(x)}dx
            \leq
                \int_{A(k,t)}
                    \left|
                        \nabla(\zeta(x)(u(x)-k)^{+})
                    \right|
                ^{p(x)}dx
        \end{equation*}

        \begin{displaymath}
            =\int_{A(k,t)}
                \left|
                    \nabla(\zeta(x)(u(x)-k))
                \right|^{p(x)}dx
        \end{displaymath}
        \begin{displaymath}
            =\int_{A(k,t)}
                \left|
                    \zeta(x)\nabla(u(x)-k)+(u(x)-k)\nabla\zeta(x)
                \right|^{p(x)}dx
        \end{displaymath}
        \begin{displaymath}
        \leq
            C\left\{
                \int_{A(k,t)}
                    \left|
                        \zeta(x)
                    \right|^{p(x)}
                    \left|
                        \nabla(u(x)-k)
                    \right|^{p(x)}dx
            \right.
        \end{displaymath}
        \begin{displaymath}
            \left.
                +\int_{A(k,t)}
                    \left|
                        (u(x)-k)\nabla\zeta(x)
                    \right|^{p(x)}dx
            \right\}
        \end{displaymath}
        \begin{displaymath}
            \leq
                C\left\{
                    \int_{A(k,t)}
                    \left|
                        \nabla u(x)
                    \right|^{p(x)}dx
                        +\int_{A(k,t)}
                    \left|
                        (u(x)-k)\nabla\zeta(x)
                    \right|^{p(x)}dx
                \right\}
        \end{displaymath}
        \begin{displaymath}
            \leq C(p,q,\lambda)
                \left\{
                    \int_{A(k,t)}
                        \left|
                            u(x)
                        \right|^{p(x)}dx
                    +\int_{A(k,t)}
                        \left|
                            (u(x)-k)\nabla\zeta(x)
                        \right|^{p(x)}dx
                \right\}
            \end{displaymath}
the last inequality can be deduced from (\ref{def:weak solution}) in which we take
    $h(x)=u(x)$
such that
    \begin{displaymath}
        \int_{\Omega}-Div
        \left(
            p(x)
            \left|
                \nabla u(x)
            \right|^{p(x)-2}\nabla u(x)
        \right)u(x)dx
    \end{displaymath}
    \begin{displaymath}
        =\lambda\int_{\Omega}b(x)
            \left|
                u(x)
            \right|^{p(x)-2}u(x)u(x)dx
    \end{displaymath}
or
    \begin{equation}
        \int_{\Omega}
            \left\{
                p(x)
                \left|
                    \nabla u(x)
                \right|^{p(x)}+a(x)
                \left|
                    u(x)
                \right|^{p(x)}
            \right\}u(x)dx
        =\lambda\int_{\Omega}b(x)
        \left|
            u(x)
        \right|^{p(x)}dx
    \end{equation}
since
    $a(x),b(x)\in L^{\infty}(\Omega),
    1<p\leq p(x)\leq q\leq p^{\ast}<+\infty$,
we get the result above.\newline Moreover,by
    $|\nabla u(x)|
    \leq
    \frac{2}{t-s}$
we obtain
\begin{eqnarray*}
    &&\int_{A(k,s)}
        \left |
            \nabla u(x)
        \right |^{p(x)}
    dx\\
    &\leq&
    C(\lambda)\int_{A(k,t)}
        \left|
            u(x)
        \right |^{p(x)}
    dx
    +
    C\int_{A(k,t)}
        \left|
            \frac{u(x)-k}{t-s}
        \right|^{p(x)}
    dx\\
    &\leq&
    C(\lambda)
    \int_{A(k,t)}
        \left|
            u(x)-k+k
        \right |^{p(x)}
    dx
    +
    C\int_{A(k,t)}
        \left|
            \frac{u(x)-k}{t-s}
        \right |^{p(x)}
    dx\\
    &\leq&
    C(\lambda)
    \int_{A(k,t)}
        \left(
            \left|
                u(x)-k
            \right|^{p(x)}
            +
            |k|^{p(x)}
        \right)
    dx
    +
    C\int_{A(k,t)}
        \left|
            \frac{u(x)-k}{t-s}
        \right |^{p(x)}
    dx\\
    &\leq&
    C(\lambda)
    \left\{
        \int_{A(k,t)}
            |k|^{p(x)}
        dx
        +
        C\int_{A(k,t)}
            \left|
                \frac{u(x)-k}{t-s}
            \right |^{p(x)}
        dx
    \right\}
\end{eqnarray*}
applying for Younger inequality and sobolev embedding theorem, we have
\begin{eqnarray*}
    &&\int_{A(k,s)}
        \left |
            \nabla u(x)
        \right |^{p}
    dx
    \leq
    \int_{A(k,s)}
        \left(
            1
            +
            \left|
                \nabla u(x)
            \right |^{p(x)}
        \right)
    dx
\\
    &&\leq
    C(\lambda)
        \left\{
            \int_{A(k,t)}
                |k|^{p(x)}
            dx
            +
            C\int_{A(k,t)}
                \left|
                    \frac{u(x)-k}{t-s}
                \right|^{p(x)}
            dx
        \right\}
        +
        |A(k,s)|
\\
    &&\leq
    C(\lambda)
        \left\{
            \int_{A(k,t)}
                \left|
                    \frac{u(x)-k}{t-s}
                \right|^{p(x)}
            dx
            +
            (1+k^{q})|A(k,t)|
        \right\}
\\
    &&\leq
    C(\lambda)
        \left\{
            \int_{A(k,t)}
                \left(
                    1
                    +
                    \left|
                        \frac{u(x)-k}{t-s}
                    \right|^{q}
                \right)
            dx
            +
            (1+k^{q})|A(k,t)|
        \right\}
\\
    &&\leq
        C(\lambda)
        \left\{
            \int_{A(k,t)}
                \left|
                    \frac{u(x)-k}{t-s}
                \right|^{q}
            dx
            +
            (1+k^{q})|A(k,t)|
        \right\}
\end{eqnarray*}
\begin{cor}\label{collar:boundary}
Suppose
    $u(x)$
 is the solution of
    $(E)_{\lambda}$, $p(x)$
 satisfies
    $1< p \leq p(x)
    \leq
     p <  p^{*}<+\infty$
in a certain open set
    $\omega$
included in
    $\Omega$,
then for any  given spherical neighborhood
    $B_{R}$ in $\omega$
and every
    $0<\rho<R<1$, $k>0$
there stands
\begin{equation}
    \int_{A(k,\rho)}
        \left|
            \nabla u(x)
        \right|^{p(x)}
    dx
    \leq C
        \left(
            \int_{A(k,R)}
                \left|
                    \frac{u(x)-k}{R-\rho}
                \right|^{p(x)}
            dx
            +
            |A(k,R)|
        \right)
\end{equation}
and
\begin{equation}
    \int_{A(k,\rho)}
        \left|
            \nabla u(x)
        \right|^{p(x)}
    dx
    \leq C
        \left(
            \int_{A(k,R)}
                \left|
                    \frac{u(x)-k}{R-\rho}
                \right|^{q}
            dx
            +
            |A(k,R)|
        \right)
\end{equation}
where
    $C=C(p,q,\lambda)$,
    $1< p^{\circ}
    \leq
    p
    \leq
    p(x)
    \leq
    p
    <
    p^{*}
    <
    +\infty$
\end{cor}
\subsection{Local Bounded of Solution}
\begin{theo}[Local Bounded on solution]\label{theorem:Local Bounded}
the solution of
    $(E)_{\lambda}$
is local bounded, that is, for every spherical
neighborhood $B_{R}$ in
    $\omega\subset\Omega$ and $0<R<1$,
there exist a certain given positive number $k>0$ such that
    $u(x)\leq k$
where
    $p(x)$
satisfies
    $1 < p \leq p(x) \leq p <  p^{*}<+\infty$.
\end{theo}
\textbf{Proof}: For fixed
    $B_{R}\subset\subset\Omega$
and
    $R\leq 1, k\geq k_{0}>0$
let
    $\vartheta_{h}=\frac{R}{2}+\frac{R}{2^{h+1}}$,
    $\overline{\vartheta}_{h}=\frac{\vartheta_{h}+\vartheta_{h+1}}{2}$,
    $k_{h}=k(1-\frac{1}{2^{h+1}})$,$h=0,1,2,\cdots$.
Obviously
    $\vartheta_{h}$
monotonously decreases to
    $\frac{R}{2}$
and
    $k_{h}$
monotonously increases to
    $k$
as
    $h\rightarrow+\infty$.
Define a functional $J_{h}$ as
    \begin{equation*}
    J_{h}=\int_{A(k_{h},\vartheta_{h})}|u(x)-k_{h}|^{p^{\ast}}dx
    \end{equation*}
Take a function
    $\xi(t)\in C^{1}([0,+\infty])$
such that
    $0\leq\xi(t)\leq 1, |\xi(t)|\leq C$(constant)
and
    $\xi(t)=1$
when
    $0\leq t\leq \frac{1}{2}$,
    $\xi(t)=0$
when
    $t\geq \frac{3}{4}$.
from above making the cutting function
    $\xi_{h}(x)$
as
    $\xi_{h}(x)=\xi(\frac{2^{h+1}}{R}(|x|-\frac{R}{2}))$
by the Poincar\'{e} inequality(\ref{ineq:pocar}) and Caccioppoli
inequality(\ref{ineq:cacci}) we have
\begin{eqnarray*}
      J_{h}
      &=&
      \int_{A(k_{h},\vartheta_{h})}
        |u(x)-k_{h}|^{p^{\ast}}
      dx
\\
    &\leq&
        \int_{A(k_{h},\overline{\vartheta}_{h})}
            |(u(x)-k_{h})\xi_{h}(x)|^{p^{\ast}}
        dx
\\
    &=&
    \int_{A(k_{h},\overline{\vartheta}_{h})}
        |(u(x)-k_{h})^{+}\xi_{h}(x)|^{p^{\ast}}
    dx
\\
    &\leq&
    \left \{
        \prod^{n}_{i=1}
        \left (
            \int_{A(k_{h},\overline{\vartheta}_{h})}
                \left|
                    \frac{\partial((u(x)-k_{h+1})^{+}\xi_{h}(x))}{\partial x_{i}}
                \right|^{p}
            dx
        \right )^{\frac{1}{p}}
    \right \}^{\frac{p^{\ast}}{n}}
\\
    &\leq&
    C\left \{
        \prod^{n}_{i=1}
            \left (
                \int_{A(k_{h},\overline{\vartheta}_{h})}
                    \left|
                        \frac{\partial u(x)}{\partial x_{i}}
                    \right|^{p}
                dx
                +
                2^{hp}
                \int_{A(k_{h},\overline{\vartheta}_{h})}
                    \left|
                        u(x)-k_{h+1}
                    \right|^{p}
                dx
            \right )^{\frac{1}{p}}
    \right \}^{\frac{p^{\ast}}{n}}
\end{eqnarray*}
\begin{eqnarray*}
    &\leq&
        C(\lambda)
        \left \{
            \prod^{n}_{i=1}
            \left (
                \int_{A(k_{h},\overline{\vartheta}_{h})}
                    \left|
                        \frac{u(x)-k_{h+1}}
                        {\vartheta_{h}-\overline{\vartheta}_{h}}
                    \right|^{p^{\ast}}
                dx
            \right.
        \right.
\\
    & &
    \left.
        \left.
            +2^{hp}
            \int_{A(k_{h},\overline{\vartheta}_{h})}
                \left|
                    u(x)-k_{h+1}
                \right|^{p}
            dx
            +
            (1+k^{p^{\ast}}_{h+1})
            \left|
                A(k_{h+1},\vartheta_{h})
            \right|
        \right )^{\frac{1}{p}}
    \right\}^{\frac{p^{\ast}}{n}}
\\
    &\leq&
    C(\lambda)
    \left\{
        \prod^{n}_{i=1}
        \left (
            2^{hp^{\ast}}
            \int_{A(k_{h},\overline{\vartheta}_{h})}
                \left|
                    u(x)-k_{h+1}
                \right|^{p^{\ast}}
            dx
            +
            (1+k^{p^{\ast}}_{h+1})
            \left|
                A(k_{h+1},\vartheta_{h})
            \right|
        \right.
    \right.
\\
    & &
    \left.
        \left.
            +2^{hp}
            \int_{A(k_{h},\overline{\vartheta}_{h})}
                \left|
                    u(x)-k_{h+1}
                \right|^{p}
            dx
        \right )^{\frac{1}{p}}
    \right \}^{\frac{p^{\ast}}{n}}
\end{eqnarray*}
or
      \begin{equation}\label{ineq:J1}
      J_{h+1}\leq C(\lambda)\left \{
        \prod^{n}_{i=1}
        \left (
            2^{hp^{\ast}}J_{h}
            +(1+k^{p^{\ast}}_{h+1})
            \left|
                A(k_{h+1},\vartheta_{h})
            \right|
            +2^{hp}
            \left|
                A(k_{h+1},\vartheta_{h})
            \right|^{\frac{p^{\ast}-p}{p}}
            J_{h}^{\frac{p}{p^{\ast}}}
        \right )^{\frac{1}{p}}
      \right \}^{\frac{p^{\ast}}{n}}
    \end{equation}
but
\begin{eqnarray*}
    &&(k_{h+1}-k_{h})^{p^{\ast}}
    \left|
        A(k_{h+1},\vartheta_{h})
    \right|
    \\
    &=&
    \int_{A(k_{h},\vartheta_{h})}
            \left|
                k_{h+1}-k_{h}
            \right|^{p^{\ast}}dx
    \leq \int_{A(k_{h},\vartheta_{h})}
            \left|
                u(x)-k_{h}
            \right|^{p^{ast}}dx
    \leq J_{h}
\end{eqnarray*}
or
\begin{eqnarray}\label{ineq:J2}
       \left|
        A(k_{h+1},\vartheta_{h})
    \right|
    \leq
    (\frac{2^{h+2}}{k})^{p^{\ast}}J_{h}
    ,
    k^{p^{\ast}}\left|
        A(k_{h+1},\vartheta_{h})
    \right|
    \leq
    (2^{h+2})^{p^{\ast}}J_{h}
\end{eqnarray}
combining(\ref{ineq:J1}) and (\ref{ineq:J2}) we get
\begin{eqnarray*}
  J_{h+1} &\leq&
            C(\lambda)
            \left\{
            \prod^{n}_{i=1}
            \left(
            2^{hp{\ast}}J_{h}
            +
            (1+\frac{1}{k^{p^{\ast}}})
            (2^{h+2})^{p{\ast}}J_{h}
            \right.
            \right.
   \\
   &&
   +
   \left.
   \left.
   \left\{
    (\frac{2^{h+2}}{k})
   \right\}^{\frac{p^{\ast}-p}{p^{\ast}}}
   J_{h}^{\frac{p^{\ast}-p}{p^{\ast}}}
   J_{h}^{\frac{p}{p^{\ast}}}
    \right)^{\frac{1}{p}}
    \right\}^{\frac{p^{\ast}}{n}}
   \\
  &=&
  C(\lambda)
  \left\{
    2^{hp{\ast}}J_{h}
    +
    (1+\frac{1}{k^{p^{\ast}}})
    2^{hp{\ast}+2p{\ast}}J_{h}
    +2^{(h+2)(p^{\ast}-p)}
    \left(
    \frac{1}{k}
    \right)^{p^{\ast}-p}
    J_{h}
  \right\}^{\frac{p^{\ast}}{n}} \\
  &\leq&
   C(\lambda , k)
   \left\{
        2^{hp{\ast}}J_{h}
        +
        2^{hp{\ast}+2p^{\ast}}J_{h}
        +
        2^{h(p{\ast}-p)}J_{h}
   \right\}^{\frac{p^{\ast}}{n}}
    \\
  &\leq&
   C(\lambda , k)
   2^{\frac{hp{\ast}^{2}}{p}}J_{h}^{\frac{p{\ast}-p}{p+1}} \\
  &=&
  C(\lambda , k)
  \left(
  2^{\frac{p{\ast}^{2}}{p}}
  \right)^{h}
  J_{h}^{\frac{p{\ast}-p}{p+1}}
\end{eqnarray*}
by choosing $k>1$ such that
\begin{equation}\label{ineq:J0}
    J_{0}
    =
    \int_{A(\frac{k}{2},R)}
    \left|
    u(x)-\frac{k}{2}
    \right|^{p^{\ast}}
    dx
    \leq
    C(\lambda)^{\frac{-1}{\eta}}
    (2^{I})^{\frac{-1}{\eta^{2}}}
\end{equation}
where
    $
    I=\frac{(p^{\ast})^{2}}{p},
    \eta=\frac{(p^{\ast})^{2}-p}{p}=\frac{p}{n-p}
    $.
Then with the help of Moser iterative inequality Lemma \ref{ineq:Moser}, we have
\begin{equation}\label{LimJ:prf of local b}
    \lim_{h\rightarrow +\infty}J_{h}=0
\end{equation}
or
\begin{equation}\label{mean value integral}
    \int_{A(k,\frac{R}{2})}
    |u(x)-k|^{p^{\ast}}
    dx
    =0
\end{equation}
This shows that
    $u(x)\leq k, x\in B_{\frac{R}{2}}$.
this completes the proof of theorem \ref{theorem:Local Bounded}

Consequently, by the similar argument to
    $-u(x)$
we can prove that
    $u(x)$
is local bounded.
    then with the compact property of
    $u(x)$
we get theorem
    \ref{theorem:Bounded of solution}
as following:
\begin{theo}[Bounded on solution]\label{theorem:Bounded of solution}
the solution of
    $(E)_{\lambda}$
is bounded.
\end{theo}
\section{H\"{o}lder Continuity of Solution}
\subsection{Sobolev-Poincar\'{e} Inequality on Solution}
\begin{lem}
\label{lemma:Sobolev-poincare inequality on solution} if
    $p(x)\in P(\Omega)$
satisfies
    $\partial p/\partial x_{i}$
is bounded almost everywhere in
    $\Omega$,
    $\sup u(x)\leq M,x\in B_{R}$,
then there exist
    $R_{1}=R_{1}(M,n,p(x))$
such that for every spherical neighborhood
    $B_{R}\subset\subset\Omega, 0<R<R_{1}$
and
    $1 < p \leq p(x) \leq p < p^{*}<+\infty, (x\in B_{R})$
the solution of eigenvalue problem
    $(E)_{\lambda}$ $u(x)$
satisfies
    \begin{equation}
    \left (\oint_{B_{R}} \left | \frac{u(x)}{R}\right
    |^{\frac{p(x)n}{n-1}}dx \right )^{\frac{n-1}{n}} \leq
    C(n,p,q)\oint_{B_{R}}|\nabla u(x)|^{p(x)}dx+C(n)|B_{R}|
    \end{equation}
\end{lem}
Proof: In lemma \ref{lemma:ineq:s-p with nonstandard} let
    $\gamma =1$,
notice that
    $
    \left|
    \{x\in B_{R}: |u(x)|>0\}
    \right|
    \leq
    |B_{R}|
    $
we can get the inequality easily.
\subsection{Harnack Inequality on Solution}
\begin{theo}
\label{theorem:harnack inequality on solution} if
    $p(x)\in P(\Omega)$
satisfies
    $\partial p/\partial x_{i}$
is bounded almost everywhere in
    $\Omega$, $1 < p \leq p(x) \leq p < p^{*}<+\infty$,
in
    $\omega\subset \Omega$, $\sup u(x)\leq M,x\in B_{R}$,
then there exist
    $R_{1}=R_{1}(M,n,p(x))$
such that for every spherical neighborhood
    $B_{R}\subset\omega, 0<R<R_{1}$
the solution of eigenvalue problem
    $(E)_{\lambda}$ $u(x)$
satisfies
    \begin{equation}
        \sup_{B_{\frac{R}{2}}} u(x)
        \leq
        C R^{p/q}
        \left(
            \left(
                \frac{|A(0,R)|}{R^{n}}
            \right)^{\beta}
            \oint_{A(0,R)}
                \left|u/R\right|^{p(x)}
            dx
            +
            R^{n}
        \right)^{1/q}
    \end{equation}
where
        $\beta > 0,\beta(\beta +1)=\frac{1}{n} , C=C(n,p,q,p(x),M)$
\end{theo}
Proof: Let
    $h<k, \frac{R}{2}\leq \rho<\sigma\leq R<R_{1}<1$
where
    $R_{1}$
is the same as
    Lemma\ref{lemma:Sobolev-poincare inequality on solution}.
making the cutting function on
    $B_{\rho}$
such that
    $|\nabla\zeta(x)|\leq \frac{4}{\sigma-\rho}$
. By H\"{o}lder inequality, we have
\begin{eqnarray}
\nonumber
    & &\oint_{A(k,\rho)}
    \left|
        \frac{u-k}{\rho}
    \right|^{p(x)}
    dx  \\
\nonumber
    &\leq&
    C(n,q)
    \left(
        \oint_{B_{R}}
            \left|
                \frac{(u-k)^{+}\zeta(x)}{R}
            \right|^{\frac{p(x)n}{n-1}}
        dx
    \right)^{\frac{n-1}{n}}
    \left(
        \frac{A(k,\rho)}{R^{n}}
    \right)^{\frac{1}{n}} \\
    \nonumber
& &\oint_{A(k,\rho)}
    \left|
        \frac{u-k}{\rho}
    \right|^{p(x)}
    dx \\
    \nonumber
&\leq&
    C(n,p,q)
    \left(
        \oint_{B_{R}}
            \left|
                \nabla
                    \left(
                    (u-k)^{+}\zeta(x)
                    \right)
            \right|^{p(x)}
        dx
        +
        \left|
            A
            \left(
            k,\frac{\sigma+\rho}{2}
            \right)
        \right|
    \right)
    \left(
        \frac{A(k,\rho)}{R^{n}}
    \right)^{\frac{1}{n}} \\
    \label{ineq:prf of Harnack6}
    \,
\end{eqnarray}
However,by assumed conditions for
    $x\in A
    \left(k,
    \frac{\sigma+\rho}{2}
    \right)$
we have

   \begin{equation}\label{ineq:prf of Harnack5}
    |\nabla(u(x)-k)^{+}\zeta(x)|
    \leq
    |\nabla u(x)|
    +
    4
    \left|
        \frac{u(x)-k}{\sigma-\rho}
    \right|\end{equation}
then it follows from
    (\ref{ineq:prf of Harnack5})
and
    (\ref{ineq:prf of Harnack6})
that
\begin{eqnarray}\label{ineq:prf of Harnack1}
    \nonumber
& &\oint_{A(k,\rho)}
    \left|
        \frac{u-k}{\rho}
    \right|^{p(x)}
    dx \\
    \nonumber
&\leq&
    C(n,p,q)
    \left(
        \frac{|A(k,\rho)|}{R^{n}}
    \right)^{\frac{1}{n}}
    \left\{
        \oint_{A(k,\frac{\sigma+\rho}{2})}
        \left(
            \left|
                \nabla u(x)
            \right|^{p(x)}
            +
            \left|
            \frac{u-k}{\sigma-\rho}
            \right|^{p(x)}
        \right)
        dx
    \right. \\
        & &
    \left.
        +
        \left|
            A
            \left(
            k,\frac{\sigma+\rho}{2}
            \right)
        \right|
    \right\}
    \end{eqnarray}
combining (\ref{ineq:prf of Harnack1}) and Lemma \ref{collar:boundary} we have
\begin{eqnarray}\label{ineq:prf of harnack2}
\nonumber
     & &\oint_{A(k,\rho)}
    \left|
        \frac{u-k}{\rho}
    \right|^{p(x)}
    dx \\
    \nonumber
    &\leq&
    C(n,p,q)
    \left(
        \frac{|A(k,\rho)|}{R^{n}}
    \right)^{\frac{1}{n}}
    \left\{
        \oint_{A(k,\sigma)}
            \left|
            \frac{u-k}{\sigma-\rho}
            \right|^{p(x)}
            dx
            +
        \left|
            A
            \left(
            k,\sigma
            \right)
        \right|
    \right\} \\
\,
\end{eqnarray}
moreover when
    $h<k$
\begin{eqnarray}
\label{ineq:prf of Harnack3}
  |A(k,\rho)|
  &\leq&
      \frac{\sigma^{p}}{(k-h)^{q}}
    \int_{A(h,\sigma)}
  \left|
    \frac{u-h}{\sigma}
  \right|^{p(x)}dx  \\
  \label{ineq:prf of Harnack4}
  \oint_{A(k,\sigma)}
  \left|
    \frac{u-k}{\sigma}
  \right|^{p(x)}dx
   &\leq&
   \oint_{A(h,\sigma)}
  \left|
    \frac{u-k}{\sigma}
  \right|^{p(x)}
  dx
\end{eqnarray}
therefore, from
    (\ref{ineq:prf of harnack2})
and
    (\ref{ineq:prf of Harnack4})
\begin{eqnarray}\label{ineq:prf of harnack8}
\nonumber
  & &\oint_{A(k,\rho)}
    \left|
        \frac{u-k}{\rho}
    \right|^{p(x)}
    dx   \\
    \nonumber
  &\leq&
  C(n,p,q)
    \left(
        \frac{|A(h,\rho)|}{R^{n}}
    \right)^{\frac{1}{n}}
    \left\{
        \oint_{A(h,\sigma)}
            \left|
            \frac{u-h}{\sigma-\rho}
            \right|^{p(x)}
            dx
            \right.\\
& &
            \left.
            +
            \frac{\sigma^{p}}{(k-h)^{h}}
        \oint_{A(h,\sigma)}
            \left|
            \frac{u-h}{\sigma-\rho}
            \right|^{p(x)}
            dx
    \right\}
\end{eqnarray}
noticing that
\begin{equation}\label{ineq:prf of harnack7}
    \left(
        \frac{|A(k,\rho)|}{R^{n}}
    \right)^{\beta}
    \leq
    \left(
        C(n)\frac{\sigma^{p}}{(k-h)^{q}}
        \oint_{A(h,\sigma)}
        \left|
        \frac{u-h}{\sigma}
        \right|^{p(x)}
        dx
    \right)^{\beta}
\end{equation}
if we multiply the two sides of (\ref{ineq:prf of harnack7}) by
    \begin{displaymath}
        \oint_{A(k,\rho)}
    \left|
        \frac{u-k}{\rho}
    \right|^{p(x)}
    dx
    \end{displaymath}
and
    \begin{displaymath}
        C(n,p,q)
    \left(
        \frac{|A(h,\rho)|}{R^{n}}
    \right)^{\frac{1}{n}}
    \left\{
        \oint_{A(h,\sigma)}
            \left|
            \frac{u-h}{\sigma-\rho}
            \right|^{p(x)}
            dx
                        +
            \frac{\sigma^{p}}{(k-h)^{h}}
        \oint_{A(h,\sigma)}
            \left|
            \frac{u-h}{\sigma-\rho}
            \right|^{p(x)}
            dx
    \right\}
    \end{displaymath}
we get
    \begin{eqnarray}\label{ineq:prf of harnack9}
    \nonumber
        & &
        \left(
            \frac{|A(k,\rho)|}{R^{n}}
        \right)^{\beta}
        \oint_{A(k,\rho)}
            \left|
                \frac{u-k}{\rho}
            \right|^{p(x)}
        dx
    \\
    \nonumber
        &\leq&
        C(n,p,q)
        \left(
            \frac{|A(h,\rho)|}{R^{n}}
        \right)^{\frac{1}{n}}
        \left(
            \frac{\sigma^{p}}{(k-h)^{q}}
            \oint_{A(h,\sigma)}
                \left|
                    \frac{u-h}{\sigma}
                \right|^{p(x)}
            dx
        \right)^{\beta}
    \\
        & &
        \left\{
            \oint_{A(h,\sigma)}
                \left|
                    \frac{u-h}{\sigma-\rho}
                \right|^{p(x)}
            dx
            +
            \frac{\sigma^{p}}{(k-h)^{h}}
            \oint_{A(h,\sigma)}
                \left|
                    \frac{u-h}{\sigma-\rho}
                \right|^{p(x)}
            dx
        \right\}
    \end{eqnarray}
Now we take
    $\sigma
    =R_{i}
    =\frac{R}{2}+\frac{R}{2^{i+1}},
    \rho
    =R_{i+1}
    ,
    h
    =k_{i}
    =dR^{\frac{p}{q}}(1-\frac{1}{2^{i}}),
    k=k_{i+1}
    $
for every
    $i\in \mathbb{N}$
and some
    $d\in \mathbb{R}$
to be chosen later. Taking into account that
    $k_{i+1}-k_{i}
    =\frac
    {
    dR^{
        \frac{p}
        {q}}
    }
    {2^{I+1}}$,
    $R_{i}-R_{i+1}
    =\frac{R}{2^{i+1}}$
from
    (\ref{ineq:prf of harnack9})
we have
\begin{eqnarray*}
   & &
   \left(
        \frac
        {|A(k_{i+1}, R_{i+1})|}
        {R^{n}}
   \right)^{\beta}
    \oint_{A(k_{i+1}, R_{i+1})}
        \left|
            \frac
            {u-k_{i+1}}
            {R_{i+1}}
        \right|^{p(x)}
    dx
   \\
   &\leq&
   C
   \frac
   {2^{(1+\beta)qi}}
   {d^{q\beta}}
    \left(
        1+\frac
        {R^{n}}
        {d^{q}}
    \right)
    \left(
        \frac
        {|A(k_{i}, R_{i})|}
        {R^{n}}
   \right)^{\frac{1}{n}}
    \left(
        \oint_{A(k_{i}, R_{i})}
        \left|
            \frac
            {u-k_{i}}
            {R_{i}}
        \right|^{p(x)}
    dx
    \right)^{1+\beta}
\end{eqnarray*}
where
    $C=C(n,p,q)$.
let
    \begin{eqnarray*}
      \varphi(k, \rho)
      &=&
        \left(
        \frac
        {|A(k, \rho)|}
        {R^{n}}
   \right)^{\beta}
    \oint_{A(k, \rho)}
        \left|
            \frac
            {u-k}
            {\rho}
        \right|^{p(x)}
    dx
    \end{eqnarray*}
we have
    \begin{eqnarray}\label{ineq:prf of harnack10}
        \varphi(k_{i+1}, R_{i+1})
        &\leq&
      C(n,p,q)
   \frac
   {2^{(1+\beta)qi}}
   {d^{q\beta}}
    \left(
        1+\frac
        {R^{n}}
        {d^{q}}
    \right)
    \varphi^{1+\beta}(k_{i}, R_{i})
    \end{eqnarray}
choosing
    $d$
such that
    $d^{q}\leq R^{n}$
and
    $\varphi(k_{0}, R_{0})
    =\varphi(0, R)
    \leq
    C(n,p,q)2^{-\frac{(1+\beta)q}{\beta^{2}}}d^{q}$
from (\ref{ineq:prf of harnack10}) and Lemma\ref{ineq:Moser}, we get
\begin{equation}\label{eq:prf of Harnack1}
    \lim_{i\rightarrow+\infty}
    \varphi(k_{i}, R_{i})
    =\varphi(dR^{\frac{p}{q}},\frac{R}{2})
    =0
\end{equation}
taking
    $d=R^{n}+C(n,p,q)\varphi(0, R)$
we deduce the desired result.
\begin{theo}\label{theorem:Holder continuity}
the weak solution of $(E)_{\lambda}$ is local H\"{o}lder continuous.
\end{theo}
Proof: by the sane proof with $-u(x)$ and notice that
    $\oint_{B_{R}}|u(x)|^{p(x)}dx$
is bounded with
    $p<p(x)< q$
we can get the estimate of $u(x)$ on $B_{R}$
    \begin{equation}\label{ineq:Holder}
        osc(u(x), \frac{R}{2})
        \leq
        C(n,p,q,M)
        R^{\frac{p}{q}}
    \end{equation}
where
    $M=\sup_{B_{R}}u(x)$.
\section{On the First Eigenvalue}
    \subsection{Comparison Principle}
        \begin{lem}[Comparison Principle]\label{Lemma:comparison principle}
            Let
                $F(x,u):\Omega \times R^{1}\rightarrow R^{1}$
            be measurable in $x$ and monotone nondecreasing in $u$, let
                $u_{1},u_{2}\in W^{1,p(x)}(\Omega)$
            satisfies
                \begin{equation}\label{ineq:compar}
                    -\Delta_{p(x)}u_{1}+F(x,u_{1})
                    \leq -\Delta_{p(x)}u_{2}+F(x,u_{2})
                \end{equation}
            in
                $W^{-1,p'(x)}(\Omega),p'(x)=p(x)/p(x)-1$
            Then
                $u_{1}\leq u_{2}$ on $\partial\Omega$
            implies
                $u_{1}\leq u_{2}$ in $\Omega$.
        \end{lem}
\textbf{Proof}:Put
    $\omega(x)=max(u_{1}-u_{2},0)$
by
    $u_{1},u_{2}\in W^{1,p(x)}(\Omega)$,
    $\omega(x)\in W^{1,p(x)}(\Omega)$.
Multiplying(\ref{ineq:compar}) by $\omega$ and using the
monotonicity of $F(x,u)$, we have
    \begin{eqnarray*}
      & &
      \int_{\Omega}
      \left(
        -\Delta_{p(x)}u_{1}(x)+F(x, u_{1}(x))
      \right)\omega(x)
      dx
      \\
       &\leq&
        \int_{\Omega}
      \left(
        -\Delta_{p(x)}u_{2}(x)+F(x, u_{2}(x))
      \right)\omega(x)
      dx
       \\
       & &
      \int_{\Omega}
      -\Delta_{p(x)}u_{1}(x)
        \omega(x)
        dx
        +
        \int_{\Omega}
        F(x, u_{1}(x))
        \omega(x)
      dx
      \\
       & \leq &
        \int_{\Omega}
        -\Delta_{p(x)}u_{2}(x)
        \omega(x)
        dx
        +
        \int_{\Omega}
        F(x, u_{2}(x))
        \omega(x)
      dx
    \end{eqnarray*}
or
    \begin{eqnarray*}
       & &
      \int_{
      \{
      x\in\Omega:u_{1}(x)\geq u_{2}(x)
      \}
      }
      -\Delta_{p(x)}u_{1}(x)
        \omega(x)
        dx
        +
        \int_{
        \{
      x\in\Omega:u_{1}(x)\geq u_{2}(x)
      \}
        }
        F(x, u_{1}(x))
        \omega(x)
      dx
      \\
       & \leq &
        \int_{
        \{
      x\in\Omega:u_{1}(x)\geq u_{2}(x)
      \}
        }
        -\Delta_{p(x)}u_{2}(x)
        \omega(x)
        dx
        +
        \int_{
        \{
      x\in\Omega:u_{1}(x)\geq u_{2}(x)
      \}
        }
        F(x, u_{2}(x))
        \omega(x)
      dx
    \end{eqnarray*}
        \begin{eqnarray*}
       & &
      \int_{
      \{
      x\in\Omega:u_{1}(x)\geq u_{2}(x)
      \}
      }
      -\Delta_{p(x)}u_{1}(x)
        \omega(x)
        dx
      \\
       & \leq &
        \int_{
        \{
      x\in\Omega:u_{1}(x)\geq u_{2}(x)
      \}
        }
        -\Delta_{p(x)}u_{2}(x)
        \omega(x)
        dx
    \end{eqnarray*}
according to the definition of
    $-\Delta_{p(x)}u(x)$
it follows from above that
    \begin{equation}\label{ineq:proof of extremvalue}
    \int_{
      D
      }
      p(x)
      (
      |
      \nabla u_{1}(x)
      |^{p(x)-2}
        \nabla u_{1}(x)
-
      |
      \nabla u_{2}(x)
      |^{p(x)-2}
        \nabla u_{2}(x)
      )
      (
\nabla u_{1}(x)-\nabla u_{2}(x)
      )
        dx
        \leq 0
    \end{equation}
where
 $D=\{
      x\in\Omega:u_{1}(x)\geq u_{2}(x)
      \}$
but
\begin{equation*}
    (
      |
      \nabla u_{1}(x)
      |^{p(x)-2}
        \nabla u_{1}(x)
-
      |
      \nabla u_{2}(x)
      |^{p(x)-2}
        \nabla u_{2}(x)
      )
      (
\nabla u_{1}(x)-\nabla u_{2}(x)
      )
\geq 0
\end{equation*}
hence,
    \begin{equation*}
        \nabla u_{1}(x)
        =
        \nabla u_{1}(x)
        ,
        x\in
        \left\{
        x\in\Omega:u_{1}(x)\geq u_{2}(x)
        \right\}
    \end{equation*}
it's means
    $\nabla \omega(x)=0$
or
    $u_{1}(x)=u_{2}(x)$
when
    $x\in
        \left\{
        x\in\Omega:u_{1}(x)\geq u_{2}(x)
        \right\}$
which implies
    $u_{1}(x)\leq u_{2}(x), x\in\Omega$.
    \begin{lem}[Extremum Principle]
        If
            $u(x)\in W^{1,p(x)}_{0}(\Omega)\cap C^{1}(\overline{\Omega})$
        satisfies
            \begin{displaymath}
                \left\{
                    \begin{array}{cc}
                    -\Delta_{p(x)}u(x)+Mu^{p(x)-1}(x) \geq 0, & \textrm{in}\ \
                    W^{-1,p(x)}(\Omega), M\geq 0\\
                    u(x)>0, & x\in \Omega \\
                    u(x)=0, & x\in \partial\Omega
                    \end{array}
                \right.
            \end{displaymath}
        Then the outer normal derivative
            $\frac{\partial u}{\partial n}$
        of $u$ is strictly negative on $\partial\Omega$.
        \end{lem}
Proof: for any given
    $x_{0}\in \partial\Omega$
and a sufficiently small
    $R>0$,
There exists
    $y\in\Omega$
such that
    $B_{2R}(y)
    \subset
    \Omega$
and
    $x_{0}\in\partial B_{2R}(y)\cap\partial\Omega$
where
    $B_{\rho}(z)=:
    \left\{
        x\in R^{n}:|z-x|<\rho
    \right\}$.
Set
    \begin{equation*}
        \upsilon(x)
        =\alpha(3R-r)^{\delta}-\alpha R^{\delta},
        r=|x-y|
    \end{equation*}
for fixed $\delta$ taking sufficiently small $\alpha, R$ such that
    \begin{eqnarray}
       & & -\Delta_{p(x)}\upsilon(x)+M\upsilon^{p(x)-1}(x)\leq 0, x\in\Omega_{R}
       \label{ineq:proof of normed vector1} \\
       & & \upsilon(x)\leq u(x), x\in\partial\Omega_{R}
        \label{ineq:proof of normed vector2}
    \end{eqnarray}
where
    $\Omega_{R}=B_{2R}(y)\setminus\overline{B_{R}(y)}$.
Now we proof that
    (\ref{ineq:proof of normed vector1})
and
    (\ref{ineq:proof of normed vector2})
is valid.
    when
    $x\in\partial B_{2R}(y)$, $\upsilon(x)\leq u(x)$
    is trivial. noticing that
    $\upsilon(x)=(2^{\delta}-1)\alpha R^{\delta}, x\in\partial B_{R}(y)$,
    $u(x)>0, x\in\partial B_{R}$.
we can get
    (\ref{ineq:proof of normed vector2})
by taking sufficiently small
    $\alpha, R$.

For
    (\ref{ineq:proof of normed vector1})
according to the chain rule of differential,
    \begin{eqnarray*}
    & &
    -Div(
        p(x)|
            \nabla(
                \alpha(3R-r)^{\delta}
                -
                \alpha R^{\delta}
            )
        |^{p(x)-2}
        \nabla(
        \alpha(3R-r)^{\delta}
        -
        \alpha R^{\delta})
    )\\
    &=& (
        \alpha\delta(3R-r)^{\delta-1}
        )^{p(x)-1}
        \sum^{n}_{i=1}\frac{\partial p(x)}{\partial r}
        \cdot
        \frac{\partial r}{\partial x_{i}} \\
    & & +
        (
        \alpha\delta(3R-r)^{\delta-1}
        )^{p(x)-1}
        \ln\alpha\delta(3R-r)^{\delta-1}
        \sum^{n}_{i=1}\frac{\partial p(x)}{\partial r}
        \cdot
        \frac{\partial r}{\partial x_{i}} \\
    & & -p(x)(\alpha\delta(3R-r)^{\delta-1})^{p(x)-1}
        \frac{p(x)-1}{\alpha\delta(3R-r)^{\delta}-1}
        \alpha\delta(\delta-1)(3R-r)^{\delta-2} \\
    & & +p(x)(\alpha\delta(3R-r)^{\delta-1})^{p(x)-1}
        (\frac{n}{r}-\frac{1}{r^{3}})\\
    &\equiv& H_{1}+H_{2}+H_{3}+H_{4}
    \end{eqnarray*}
By assumed conditions,
    $\frac{\partial p(x)}{\partial x_{i}}$
and
    $\frac{\partial r}{\partial x_{i}}$ are both bounded for every
    $i=1,\ldots,n$.
hence
    $
    H_{1}\leq C_{1}(\alpha\delta(3R-r)^{\delta-1})^{p(x)-1}
    \leq C_{2}(\alpha\delta R^{\delta-1})^{p(x)-1}
    $
Taking sufficient small $R$ such that
    $\ln\alpha\delta(3R-r)^{\delta-1}\leq 0, \frac{n}{r}-\frac{1}{r^{3}}\leq 0$
then we get
    $
    H_{2}\leq 0, H_{4}\leq 0
    $.
because $p(x)$ is bounded in $\Omega$, we have
    $
    H_{3}\leq -C_{3}(\alpha\delta R^{\delta-1})^{p(x)-1}\frac{\delta-1}{3R-r}
    $.
And noticing that
    $M\upsilon (x)^{p(x)-1}=M(\alpha\delta(3R-r_{\xi})^{\delta-1}(2R-r))^{p(x)-1}$
where
    $r\leq r_{\xi}\leq 2R$
 we also have
    $M\upsilon (x)^{p(x)-1}
    \leq
    M_{1}(\alpha\delta R^{\delta-1})^{p(x)-1}
    $
 stands for sufficiently small
    $2R-r$
 All the conditions shown above imply that
\begin{eqnarray*}
  & &
  -\Delta_{p(x)}\upsilon(x)+M\upsilon^{p(x)-1}(x) \\
  &\leq&
  C_{2}(\alpha\delta R^{\delta-1})^{p(x)-1}
  -C_{3}(\alpha\delta R^{\delta-1})^{p(x)-1}\frac{\delta-1}{3R-r}
  +M_{1}(\alpha\delta R^{\delta-1})^{p(x)-1}
  \end{eqnarray*}
Let
    \begin{equation*}
    C_{2}(\alpha\delta R^{\delta-1})^{p(x)-1}
    -C_{3}(\alpha\delta R^{\delta-1})^{p(x)-1}\frac{\delta-1}{3R-r}
    +M_{1}(\alpha\delta R^{\delta-1})^{p(x)-1}
    \leq 0\end{equation*}
then we get
    \begin{equation*}
    C_{4}\frac{\delta-1}{3R-r}
    \geq
    C_{2}+M_{1}\end{equation*}
it stands for sufficiently small $R$. For fixed $\delta$ we prove that
    (\ref{ineq:proof of normed vector1})
and
    (\ref{ineq:proof of normed vector2})
are valid. according to Lemma
    \ref{Lemma:comparison principle}
we get the desired result.
\subsection{The Eigenvalue Problem of Solution}
We give some definitions about $(E)_{\lambda}$ as following
    \begin{defi}\label{def:the first eigenvalue}
    the first eigenvalue of
    $(E)_{\lambda_{1}}$
    is
    \begin{equation}
    \frac{1}{\lambda_{1}}
        =
        \sup
        \left
        \{
        R(\upsilon):=\frac{B(\upsilon)}{A(\upsilon)}
        \right |
        \left.
        \upsilon\in W=:W^{1,p(x)}_{0}(\Omega)\backslash\{0\}
        \right
        \}.
    \end{equation}
     \end{defi}
where
\begin{eqnarray*}
    &&A(\upsilon)
    =
    \int_{\Omega}
    \left\{
        |\nabla \upsilon (x)|^{p(x)}
        +
        \frac{a(x)}{p(x)}
        |\upsilon (x)|^{p(x)}
    \right\}
    dx
  \\
  &&B(\upsilon)
  =
  \int_{\Omega}\frac{b(x)}{p(x)}|\upsilon (x)|^{p(x)}dx
\end{eqnarray*}
     \begin{theo}[Boundedness of the First Eigenvalue ]\label{theorem:bound of eigen}
        For the first eigenvalue,there exists
            $C_{1},C_{2}>0$
        such that
            $C_{1}<\lambda_{1}<C_{2}$.
        \end{theo}
proof: suppose $B(u)$ is no positive for all $u\in W_{0}^{1,p(x)}(\Omega)$,
    then there exists a function sequence
    $f_{n}$
in
    $W_{0}^{1,p(x)}(\Omega)$
such that
    $f_{n}(x)\leq 0$
and
    $f_{n}(x)\rightarrow b^{+}(x)=:\max\{b(x),0\},(n\rightarrow\infty)$
which implies
    $B(b^{+}(x))\leq 0$
so
    $b^{+}\equiv 0$
it contradict with the definition of
    $b(x)$
therefore there must be exist a function
    $u_{0}(x)$
such that
    \begin{equation}\label{ineq:prf bound of eigen1}
        B(u_{0}(x))>0
    \end{equation}
or
    \begin{equation}\label{ineq:prf bound of eigen2}
        0<\lambda_{1}<\frac{1}{R(u_{0}(x))}
    \end{equation}
On the other hand, by Lemma \ref{ineq:s-p} and the definition of the first eigenvalue
we have
    \begin{equation}\label{ineq:proof boundof eigen3}
        \frac{1}{\lambda_{1}}
        \leq
        C(\Omega)\|b(x)\|_{L^{\infty}}
    \end{equation}
So we have
    \begin{equation}\label{ineq:prf bound of eige4}
        \frac{1}
        {C(\Omega)\|b(x)\|_{L^{\infty}}}
        \leq
        \lambda_{1}
        \leq
        \frac{1}
        {R(u_{0}(x))}
    \end{equation}
        \begin{theo}\label{theorem:exist of solution}
        There exists
            $u\in W_{0}^{1,p(x)}$
        so that
            $J_{\lambda_{1}}(u)=0$
        implies
            $u$
        is the solution of eigenvalue problem $(E)_{\lambda}$.
where
    \begin{equation*}
    J_{\lambda_{1}}(u)
    =
    \int_{\Omega}
        \left(
        |\nabla u(x)|^{p(x)}
        +
        \frac{a(x)}{p(x)}|u(x)|^{p(x)}
        \right)
        dx
        - \lambda_{1}
        \int_{\Omega}
        \frac{b(x)}{p(x)}|u(x)|^{p(x)}
        dx
    \end{equation*}
 \end{theo}

Proof: make the Fr\'{e}chet derivation of $J_{\lambda_{1}}(u)$ we have
    \begin{eqnarray*}
      J^{'}_{\lambda_{1}}(u)
      &=& \frac
      {
      d
       \left\{
        \int_{\Omega}
        \left(
        |\nabla u(x)|^{p(x)}
        +
        \frac{a(x)}{p(x)}
        |u(x)|^{p(x)}
        \right)
        dx
        - \lambda_{1}
        \int_{\Omega}
        \frac{b(x)}{p(x)}
        |u(x)|^{p(x)}
        dx
      \right\}
      }
      {du(x)}
     \\
      &=&
      \lim_{t\rightarrow 0}
      \frac{
        \int_{\Omega}|\nabla (u(x)+th(x))|^{p(x)}dx
        -
        \int_{\Omega}|\nabla u(x)|^{p(x)}dx
      }
        {th(x)}\\
      &&+
      \lim_{t\rightarrow 0}
      \frac{
      \int_{\Omega}
      \frac{a(x)}{p(x)}
      |u(x)+th(x)|^{p(x)}dx
      -
      \int_{\Omega}
      \frac{a(x)}{p(x)}
      |u(x)|^{p(x)}dx
      }
      {th(x)}\\
      &&-
      \lambda_{1}\lim_{t\rightarrow 0}
      \frac{
      \int_{\Omega}
      \frac{b(x)}{p(x)}
      |u(x)+th(x)|^{p(x)}dx
      -
      \int_{\Omega}
      \frac{b(x)}{p(x)}
      |u(x)|^{p(x)}dx}
      {th(x)}\\
      &=& G_{1}+G_{2}+G_{3}=0
    \end{eqnarray*}
where
    $h(x)\in C_{0}^{\infty}(\Omega)$
and
    \begin{eqnarray}
    \nonumber
      G_{1} &=&
      \lim_{t\rightarrow\infty}
      \frac
      {\int_{\Omega}|\nabla(u(x)+th(x))|^{p(x)}dx
      -
      \int_{\Omega}|\nabla u(x)|^{p(x)}dx}
      {th(x)}\\
      \nonumber
      &=&
      \frac
      {\int_{\Omega}
      \lim_{t\rightarrow\infty}
      \frac{1}{t}
    \left(
      |\nabla(u(x)+th(x))|^{p(x)}
      -
      |\nabla u(x)|^{p(x)}
        \right)
      dx}
      {h(x)}\\
      \nonumber
      &=&
      \frac{1}{h(x)}
      \int_{\Omega}
      \left.
      \frac{d|\nabla(u(x)+th(x))|^{p(x)}}
      {dt}
      \right|_{t=0}
      dx\\
      \nonumber
      &=&
      \frac{1}{h(x)}
      \int_{\Omega}
      p(x)|\nabla u(x)|^{p(x)-2}
      \nabla u(x)
      \nabla h(x)
      dx\\
      \label{eq:prf of exi0}
      &=&
      \frac{1}{h(x)}
      \int_{\Omega}
        -Div
        \left(
      p(x)|\nabla u(x)|^{p(x)-2}
      \nabla u(x)
        \right)
      h(x)
      dx
    \end{eqnarray}
\begin{eqnarray}
\nonumber
    G_{2} &=&
    \lim_{t\rightarrow 0}
    \frac{
    \int_{\Omega}
    \frac{a(x)}{p(x)}|u(x)+th(x)|^{p(x)}
    dx
    -
    \int_{\Omega}
    \frac{a(x)}{p(x)}|u(x)|^{p(x)}
    dx}
    {th(x)}\\
    &=&
    \label{eq:prf of exi1}
    \frac{1}{h(x)}
    \int_{\Omega}
    a(x)|u(x)|^{p(x)-1}h(x)
    dx
\end{eqnarray}
\begin{eqnarray}
\nonumber
    G_{3}
    &=&
    \lim_{t\rightarrow 0}
    \frac{
        \int_{\Omega}
            \frac{b(x)}{p(x)}|u(x)+th(x)|^{p(x)}
        dx
        -
        \int_{\Omega}
            \frac{b(x)}{p(x)}|u(x)|^{p(x)}
        dx
        }
    {th(x)}\\
    &=&
    \label{eq:prf of exi2}
    \frac{1}{h(x)}
    \int_{\Omega}
    b(x)|u(x)|^{p(x)-1}h(x)dx
\end{eqnarray}
combining
    (\ref{eq:prf of exi0}), (\ref{eq:prf of exi1}) and (\ref{eq:prf of exi2})
we get
    \begin{eqnarray*}
    &&
    \int_{\Omega}
        -Div\left(
            p(x)|\nabla u(x)|^{p(x)-2}
            \nabla u(x)
        \right)
        h(x)
    dx
    +
    \int_{\Omega}
        a(x)|u(x)|^{p(x)-1}h(x)
    dx
    \\
    &&
    =
    \lambda_{1}
    \int_{\Omega}
        b(x)|u(x)|^{p(x)-1}h(x)
    dx
    \end{eqnarray*}
therefore,
    $u(x)$
is the solution of eigenvalue problem
    $(E)_{\lambda}$.
Noticing that
    $J_{\lambda_{1}}(u(x))
    =
    J_{\lambda_{1}}(|u(x)|)
    $
and
    $W^{1,p(x)}_{0}(\Omega)
    \hookrightarrow
    L^{p(x)}(\Omega)$
is a continuously compact imbedding. we obtain the existence of nonegative solution
about the eigenvalue problem $(E)_{\lambda_{1}}$.
\begin{theo}\label{theorem:simple of solv}
    The first eigenvalue
    $\lambda_{1}$
    is simple, that is to say, the set of solutions is
    $\{tu(x):t\in$$\mathbb{R}\}$
\end{theo}
Proof: the sufficient is obvious. Let
    $u_{1}, u_{2}$
be two solutions of
    $(E)_{\lambda_{1}}$
and
    $M(t,x)=\max(u_{1},tu_{2})$, $m(t,x)=\min(u_{1}, tu_{2})$.
because
\begin{equation*}
    \frac{1}{\lambda_{1}}=\sup
    \{
    R(\upsilon):=\frac{B(\upsilon)}{A(\upsilon)};
    \upsilon\in W:=W_{0}^{1,p(x)}(\Omega)\backslash\{0\}
    \}
    \geq
    \frac{B(M)}{A(M)}
\end{equation*}
therefore
    \begin{equation}\label{ineq:prf of simp1}
    J_{\lambda_{1}}(M)=A(M)-\lambda_{1}B(M)\geq 0
    \end{equation}
In the same like,
    \begin{equation}\label{ineq:prf of simp2}
    J_{\lambda_{1}}(m)=A(m)-\lambda_{1}B(m)\geq 0
    \end{equation}
 Now we show that
    $M,m$
are the solutions of
    $(E)_{\lambda_{1}}$
or
    $J_{\lambda_{1}}(M)
    =
    J_{\lambda_{1}}(m)=0$.
By the definition of
    $J_{\lambda_{1}}$
it's easy to see that
\begin{eqnarray*}
    &&
    J_{\lambda_{1}}(M)
    +
    J_{\lambda_{1}}(m)\\
    &=&
    A(M)-\lambda_{1}B(M)
    +
    A(m)-\lambda_{1}B(m)\\
    &=&
     \int_{\Omega}
        \left(
        |\nabla M|^{p(x)}+a|M|^{p(x)}
        \right)
        dx
        - \lambda_{1}
        \int_{\Omega}
        b|M|^{p(x)}
        dx\\
        &&
        +
        \int_{\Omega}
        \left(
        |\nabla m|^{p(x)}+a|m|^{p(x)}
        \right)
        dx
        -\lambda_{1}
        \int_{\Omega}
        b|m|^{p(x)}
        dx\\
        &=&
     \int_{I_{u_{1}}}
        \left(
        |\nabla M|^{p(x)}+a|M|^{p(x)}
        \right)
        dx
        - \lambda_{1}
        \int_{I_{u_{1}}}
        b|M|^{p(x)}
        dx\\
        &&
        +
        \int_{I_{u_{1}}}
        \left(
        |\nabla m|^{p(x)}+a|m|^{p(x)}
        \right)
        dx
        -\lambda_{1}
        \int_{I_{u_{1}}}
        b|m|^{p(x)}
        dx\\
        &&
        +
     \int_{I_{tu_{2}}}
        \left(
        |\nabla M|^{p(x)}+a|M|^{p(x)}
        \right)
        dx
        - \lambda_{1}
        \int_{I_{tu_{2}}}
        b|M|^{p(x)}
        dx\\
        &&
        +
        \int_{I_{tu_{2}}}
        \left(
        |\nabla m|^{p(x)}+a|m|^{p(x)}
        \right)
        dx
        -\lambda_{1}
        \int_{I_{tu_{2}}}
        b|m|^{p(x)}
        dx\\
        &=&
     \int_{\Omega}
        \left(
        |\nabla u_{1}|^{p(x)}+a|u_{1}|^{p(x)}
        \right)
        dx
        - \lambda_{1}
        \int_{\Omega}
        b|u_{1}|^{p(x)}
        dx\\
        &&
        +
        \int_{\Omega}
        \left(
        |\nabla u_{2}|^{p(x)}+a|u_{2}|^{p(x)}
        \right)
        dx
        -\lambda_{1}
        \int_{\Omega}
        b|u_{2}|^{p(x)}
        dx\\
        &=&
        J_{\lambda_{1}}(u_{1})
        +J_{\lambda_{1}}(u_{2})
        =0
        \end{eqnarray*}
where
    $I_{u_{1}}=\{x\in\Omega: u_{1}(x)\geq u_{2}(x)\}$,
    $I_{tu_{2}}=\{x\in\Omega: u_{1}(x)<u_{2}(x)\}$
according to
    (\ref{ineq:prf of simp1})
and
    (\ref{ineq:prf of simp2})
    $u_{1}, tu_{2}$
be two solutions of
    $(E)_{\lambda_{1}}$.
From theorem
    \ref{theorem:Holder continuity}
we have
    $M\in C^{1,\theta}_{0}(\Omega)$
for all
    $t\geq 0$.

For certain given
    $x_{0}\in \Omega$,
we take
    $t_{0}=\frac{u_{1}(x_{0})}{u_{2}(x_{0})}$.
As we know, for every vector
    $e$
there stands
    \begin{eqnarray*}
    && u_{1}(x_{0}+he)-u_{1}(x_{0})
    \leq
    \max(
        u_{1}(x_{0}+he),
        t_{0}
        u_{2}(x_{0}+he)
    )
    -
    u_{1}(x_{0}) \\
    &=& M(t_{0}, x_{0}+he)-M(t_{0}, x_{0})
    \end{eqnarray*}
therefore, the partial derivative of
    $u_{1}(x)$
and
    $M(t,x)$
at
    $x_{0}$
satisfies
    \begin{eqnarray*}
      &&\frac{\partial u_{1}(x_{0})}{\partial x_{i}}
      =\lim_{h\rightarrow 0+}
      \frac{u_{1}(x_{0}+he_{i})-u_{1}(x_{0})}
      {h}\\
      &\geq&
    \lim_{h\rightarrow 0+}
      \frac{M(x_{0}+he_{i})-u_{1}(x_{0})}
      {h}
      =\frac{\partial M(x_{0})}{\partial x_{i}},
      \\
       &&\frac{\partial u_{1}(x_{0})}{\partial x_{i}}
      =\lim_{h\rightarrow 0-}
      \frac{u_{1}(x_{0}+he_{i})-u_{1}(x_{0})}
      {h}\\
      &\leq&
    \lim_{h\rightarrow 0+}
      \frac{M(x_{0}+he_{i})-u_{1}(x_{0})}
      {h}
      =\frac{\partial M(x_{0})}{\partial x_{i}}
      \\
    \end{eqnarray*}
or
\begin{equation}\label{eq:prf of simp}
    \frac{\partial u_{1}(x_{0})}{\partial x_{i}}
    =
    \frac{\partial M(x_{0})}{\partial x_{i}}
\end{equation}
 where
    $e_{i} (i=1,2,\cdots,n)$
are the unit normal vectors. Moreover, we have
    $\nabla_{x}u_{1}(x_{0})
=
     \nabla_{x}M(t_{0},x_{0})$,
where
    $\nabla_{x}=
    \left(
    \frac{\partial}{\partial x_{1}},
    \frac{\partial}{\partial x_{2}},
    \cdots ,
    \frac{\partial}{\partial x_{n}}
    \right)$
By the same way, we obtain
    $\nabla_{x}t_{0}u_{2}(x_{0})
=
     \nabla_{x}M(t_{0},x_{0})$.
Hence, the gradient of
    $\frac{u_{1}}{t_{0}u_{2}}$
at
    $x_{0}$
is
    \begin{eqnarray*}
      \nabla_{x}
      \left(
      \frac{u_{1}(x_{0})}
      {t_{0}u_{2}(x_{0})}
        \right) &=&
        \frac
        {
        u_{2}(x_{0})\nabla_{x}u_{1}(x_{0})
        -
        u_{1}(x_{0})\nabla_{x}u_{2}(x_{0})
        }
        {
        u_{1}(x_{0})-u_{2}(x_{0})
        }=0
    \end{eqnarray*}
By the arbitrary of $x_{0}$, we get
\begin{equation}\label{eq:prf of simp3}
    \frac{u_{1}(x)}{u_{2}(x)}\equiv const, x\in\Omega
\end{equation}

\begin{theo}\label{theorem:nonexi of solv}
    $(E)_{\lambda}$
has no solution for
    $\lambda >\lambda_{1}$.
\end{theo}
Proof: By theorem
    \ref{theorem:exist of solution},
we need only verify it for the positive solution.\\
Let
    $u, \upsilon$
be positive solutions of
    $(E)_{\lambda_{1}}$ and $(E)_{\lambda}$
respectively. Assume
    $b(x)\geq 0$,
from above we can select some solutions
    $u, \upsilon$
such that
    $u\leq\upsilon$
for all
    $x\in\Omega$.
then we deduce that there must exist
    $0<\eta<1$
so that
    \begin{equation}\label{ineq:Prf of nonexis0}
    -\Delta_{p(x)}u+au^{p(x)-1}
    \leq
    -\Delta_{p(x)}(\eta\upsilon)+a(\eta\upsilon)^{p(x)-1}
    \end{equation}
By the definition of solution and
    $u\leq\upsilon$
,because for all
    $h\in C^{\infty}_{0}(\Omega)$
and
    $h>0$
    \begin{eqnarray}\label{ineq:prf pf nonexi1}
        &&\int_{\Omega}
            -Div(p(x)|\nabla u|^{p(x)-2}\nabla u h
        dx
        +
        \int_{\Omega}
            a|u|^{p(x)-2}u h
        dx \nonumber \\
        &&=
        \lambda_{1}
        \int_{\Omega}
            b|u|^{p(x)-2}u h
        dx \nonumber\\
        &&\leq
        \lambda_{1}
        \int_{\Omega}
            b|\upsilon|^{p(x)-2}\upsilon h
        dx \nonumber\\
        &&=
        \frac{\lambda_{1}}{\lambda}
        \int_{\Omega}
            -Div(p(x)|\nabla \upsilon|^{p(x)-2}\nabla \upsilon h
        dx
        +
        \frac{\lambda_{1}}{\lambda}
        \int_{\Omega}
            a|\upsilon|^{p(x)-2}\upsilon h
        dx \nonumber
    \end{eqnarray}
We need only to prove that
    \begin{eqnarray*}
    \frac{\lambda_{1}}{\lambda}
        \int_{\Omega}
            |\nabla \upsilon|^{p(x)-2}\nabla \upsilon\nabla h
        dx
        +
    \frac{\lambda_{1}}{\lambda}
        \int_{\Omega}
            \frac{a}{p(x)}|\upsilon|^{p(x)-2} h
        dx \nonumber\\
        \leq
        \int_{\Omega}
            |\nabla (\eta\upsilon)|^{p(x)-2}\nabla (\eta\upsilon)\nabla h
        dx
        +
       \int_{\Omega}
            \frac{a}{p(x)}|\eta\upsilon|^{p(x)-2} h
        dx \nonumber\\
    \end{eqnarray*}
for certain
    $0<\eta<1$.
We obtain the desired result by taking
    $\inf\eta^{p(x)-1}\geq\frac{\lambda_{1}}{\lambda}$.
therefore applying for Lemma
    \ref{Lemma:comparison principle}
    $u\leq\eta\upsilon$
in
    $\Omega$.
Repeating this prcedure, we deduce that
    $u\leq\eta^{n}\upsilon$
in
    $\Omega$
for all
    $n\in \mathbb{N}$,
which follows
    $u\equiv 0$.
This is a contradiction.

For general case, let
    $B^{+}(x)=:\max(b(x),0), b^{-}(x)=:\max(-b(x),0)$.
Then above result implies the equation
    \begin{equation*}
        -\Delta_{p(x)}\omega(x)
        +
        \left\{
            a(x)+\lambda b^{-}(x)
        \right\}
        \omega(x)^{p(x)-1}
        =
        \mu b^{+}(x)\omega^{p(x)-1}
    \end{equation*}
has a nontrivial positive solution
    $\omega$
if
    \begin{equation*}
        \mu\leq\mu_{1}=\lambda_{1}(a(x)+\lambda b^{-}(x), b^{+}(x))
    \end{equation*}
and
    \begin{eqnarray*}
        I_{\mu_{1}}(\omega(x))
            &=&
            A(\omega(x))
            +
            \lambda
            \int_{\Omega}
                \frac{b^{-}(x)}{p(x)}
                |\omega(x)|^{p(x)}
            dx
            -
            \mu_{1}
            \int_{\Omega}
                \frac{b^{+}(x)}{p(x)}
                |\omega(x)|^{p(x)}
            dx \\
            &=&
            \min\left\{
            I_{\mu_{1}}(z(x));z(x)\in W
            \right\}\\
            &=& 0
    \end{eqnarray*}
Since
    $\upsilon$
is a positive solution of the above equation with
    $\mu_{1}=\lambda$
we deduce that
    $\lambda\leq\mu_{1}$
and
    \begin{equation*}
        J_{\lambda}(\upsilon)
        =
        I_{\lambda}(\upsilon)
        \geq
        I_{\mu_{1}}(\upsilon)
        =
        \min\left\{
            I_{\mu_{1}}(z(x));z(x)\in W
            \right\}
        =0
    \end{equation*}
However,
    \begin{equation*}
        J_{\lambda}(u)
        =
        J_{\lambda_{1}}(u)
        -
        (\lambda-\lambda_{1})B(u)
        <0
    \end{equation*}
for all
    $u\in W^{1,p(x)}_{0}$.
This is a contradiction.

\section*{Remark}
In this chapter ,we give many properties of the solutions of
$(E)_{\lambda}$ in the sense of weak. The most important part we
discussed are the boundedness and H\"{o}lder continuity of the weak
solutions, which are also important to weak solution. Therefore,the
eigenvalue study in the partial differential operators, in essence,
includes the solution research. Finally, we show that the conditions
restrict to $p(x)$ are not necessary.it can be replaced by$1<p\leq
p(x)\leq q<+\infty$.
\chapter{Conclusion}
In this paper we introduced certain $p(x)$-Laplace operator under
the generalized Sobolev Space and proved some conclusions as
following:

We first proved the existence,boundary and H\"{o}lder continuity of
the solutions about the $p(x)$-Laplace equation which generalizes
the result of M.\^{o}tani and T.Teshima and in the same time we
obtained two inequalities: the Caccioppoli inequality and the
Harnack inequality about the equation with $p(x)$ exponent which are
two basic conclusions in the research of equations with nonstandard
exponent condition.

On the other hand, about the eigenvalue of $p(x)$-Laplace equation
we showed some properties:All eigenvalues are bounded. Solutions
about the first eigenvalue $\lambda_{1}$ is simple i.e. all the
nontrivial solutions form an one dimension subset of the space
$W^{1,p(x)}_{0}(\Omega)$. The equation has no solution with
$\lambda>\lambda_{1}$. Furthermore,we show that the outer norm
vector of the positive solutions is strictly negative at the
boundary $\partial\Omega$ and got the comparison principle.
Consequently, by the similar argument we generalized the result of
M.\^{o}tani and T.Teshima.

From the view of this point, how about the second eigenvale, the
third,$\cdots$ and the solution about them? this is our aftertimes
works.
\chapter*{Acknowlege}
In this article, my supervisor, YongQiang Fu, has given me so many advices that I can
finish it completely.I give my best grateful thanks to him. Some relative problem
about the generalized Sobolev Space has been studied in several articles written by
him.

\end{document}